\documentclass[thmsa,final,11pt]{article}

\usepackage{eurosym}
\usepackage{amsfonts}
\usepackage{amsmath}
\usepackage{amssymb}
\usepackage{mathrsfs} 
\usepackage{enumerate}
\usepackage{graphicx}
\usepackage{subfigure}
\usepackage{color,xcolor}
\usepackage[onehalfspacing,nodisplayskipstretch]{setspace}

\usepackage{xr}

\newtheorem{theorem}{Theorem}

\newtheorem{lemma}[theorem]{Lemma}

\newtheorem{proposition}[theorem]{Proposition}

\newenvironment{proof}[1][Proof]{\noindent\textbf{#1.} }{\ \rule{0.5em}{0.5em}}
\normalsize
\newenvironment{prooft}[1]{\vskip 2mm\noindent {\bf Proof of #1.}}
                    {\hfill $\square$ \vskip 2mm \noindent}

\begin{document}

\title{Limit theorems for squared increment sums of the maximum of two isotropic fractional Brownian
fields under fixed-domain asymptotics}
\author{Nicolas CHENAVIER\thanks{%
Universit\'{e} du Littoral C\^{o}te d'Opale, 50 rue F. Buisson 62228 Calais.
nicolas.chenavier@univ-littoral.fr}\quad and\quad Christian Y.\ ROBERT%
\thanks{%
1. Universit\'{e} de
Lyon, Universit\'{e} Lyon 1, Institut de Science Financi\`{e}re et
d'Assurances, 50 Avenue Tony Garnier, F-69007 Lyon, France.
2. Laboratory in Finance and Insurance - LFA CREST - Center for Research in
Economics and Statistics, ENSAE, Palaiseau, France;  christian.robert@univ-lyon1.fr%
}}
\maketitle

\abstract{We study squared increment sums of the pointwise maximum of two independent
and identically distributed isotropic fractional Brownian fields over a fixed
two-dimensional domain. The fields are observed at the points of a homogeneous
Poisson point process with intensity \(N\), and increments are computed along
the edges of the associated Delaunay triangulation. In contrast with the case
of a single fractional Brownian field, where centered squared increment sums
satisfy a central limit theorem after the usual normalization, the pointwise
maximum exhibits a different asymptotic regime. The dominant contribution
comes from Delaunay edges located in a shrinking neighborhood of the random
interface where the two fractional Brownian fields exchange the role of the
maximizer. For Hurst parameter \(H<1/2\), we prove that the properly
normalized squared increment sum converges in probability to a deterministic
constant times the local time at zero of the difference between the two fields.
This shows that the asymptotic behavior is governed by the geometry of the
random contact set rather than by Gaussian fluctuation effects. The result
provides a key ingredient for fixed-domain asymptotic inference in
Brown--Resnick type models based on randomly located observations.}


\strut

\textit{AMS (2020):} 60F05, 60G22, 60D05, 60G55.

\maketitle

\section{Introduction}

Let \(W\) be an isotropic fractional Brownian field on \(\mathbb R^2\), with
Hurst parameter \(H<1/2\), observed on a random set of sites in the fixed
square \(\mathbf C=(-1/2,1/2]^2\). More precisely, the underlying observation sites are generated by a homogeneous Poisson point process with intensity \(N\) on \(\mathbb R^2\), and the quadratic variation is computed over Delaunay edges oriented from their leftmost endpoint, retaining only those edges whose initial endpoint belongs to the fixed window \(\mathbf C\). The asymptotic regime considered throughout the paper is the fixed-domain, or infill, regime \(N\to\infty\). In \cite{Chenavier&Robert25a}, this framework was used to establish central
limit theorems for centered squared increment sums of a single isotropic
fractional Brownian field. In particular, when the increments are computed
along the edges of the Delaunay triangulation generated by the Poisson point
process, the corresponding normalized quadratic variation satisfies a
Gaussian limit theorem. The convergence rates obtained there are consistent
with those arising for increment-based statistics on regular grids, as in
Theorem 3.2 of \cite{Chan&Wood00} and Theorem 1 of \cite{Zhu&Stein02}.

The present paper shows that this Gaussian fluctuation picture changes
substantially when the observed field is no longer a single fractional
Brownian field but the pointwise maximum of two independent copies. More
precisely, let \(W^{(1)}\) and \(W^{(2)}\) be two independent and identically
distributed isotropic fractional Brownian fields, and define
\[
        W_\vee(x)=W^{(1)}(x)\vee W^{(2)}(x),
        \qquad x\in\mathbb R^2 .
\]
We study the centered squared increment sum of \(W_\vee\) along the Delaunay
edges. Although this statistic is formally close to the one considered in
\cite{Chenavier&Robert25a}, its asymptotic behavior is of a different nature.
The reason is that, at small scales, the maximum field locally behaves like
one of the two underlying fractional Brownian fields only away from the
random interface
\[
        \{x\in\mathbf C: W^{(1)}(x)=W^{(2)}(x)\}.
\]
On each side of this interface, the local behavior of \(W_\vee\) is that of
one of the two underlying fractional Brownian fields. Near the interface,
however, the identity of the maximizer may change between the two endpoints
of a short Delaunay edge, producing transition terms which dominate the
asymptotics.

Our main result identifies this transition mechanism precisely. We prove that,
after a normalization which differs from the one in the Gaussian central limit
theorem for a single fractional Brownian field, the centered squared increment
sum of \(W_\vee\) converges in probability to a non null constant times
the local time at zero of the difference field
\[
        W^{(2\backslash 1)}(x)=W^{(2)}(x)-W^{(1)}(x).
\]
Thus the limiting object is not a centered Gaussian random variable, but a
geometric functional measuring the size, in the occupation-density sense, of
the random contact set where the two fields coincide. This result reveals a
new asymptotic regime for quadratic variation statistics of non-Gaussian
fields obtained through pointwise maxima: the leading contribution is not
produced by the usual accumulation of dependent Gaussian increment
fluctuations, but by the edges located in a shrinking neighborhood of the
interface between the two competing fields.

The use of the Delaunay triangulation is motivated by the statistical problem underlying this work. In fixed-domain inference for Brown--Resnick max-stable random fields, composite likelihood estimators can be built from pairs and triples of nearby observation sites selected through the Delaunay triangulation (see \cite{Chenavier&Robert25c}). This construction provides a data-driven notion of local neighborhood for irregularly located sites and avoids imposing an artificial regular grid. It is also geometrically natural, since the Delaunay triangulation is one of the most regular triangulations associated with a given point configuration, in particular through its angle-optimality property. In the Brown--Resnick setting, the relevant local likelihood expansions involve small spatial increments of Gaussian spectral processes, and the analysis of pairwise and triplewise contributions leads naturally to squared increment sums along Delaunay edges. For a single fractional Brownian field, this gives rise to Gaussian fluctuation limits, as studied in \cite{Chenavier&Robert25a}. The present paper addresses the next step needed for max-stable inference: the behavior of such Delaunay-edge statistics when the observed field is the pointwise maximum of two independent fractional Brownian fields. The Poisson assumption provides a tractable infill model for irregular observation sites and allows explicit expectation and covariance computations through the Slivnyak--Mecke formula.

Compared with \cite{Chenavier&Robert25a}, the contribution of the present
paper is therefore threefold. First, we show that the normalization of the
Delaunay-edge squared increment sum must be changed when passing from a single
fractional Brownian field to the pointwise maximum of two such fields.
Second, we prove that the limiting behavior is governed by the local time at
zero of \(W^{(2\backslash 1)}\), which quantifies the amount of contact
between the two fields inside the observation window. Third, we separate the
quadratic variation into two components: a Gaussian-type contribution coming
from regions where the same field remains the maximizer at both endpoints of
an edge, and a transition contribution coming from edges whose starting point
lies in a shrinking neighborhood of the interface where the two fields
exchange the role of the maximizer. We show that, at the relevant scale, the
first component is negligible while the second one converges to the
local-time limit. This decomposition is the central mechanism behind the
non-Gaussian asymptotic regime obtained in this paper.

The proof strategy combines tools from stochastic geometry and Gaussian
analysis. The Poisson--Delaunay structure is used to control the distribution
of edge lengths and to replace local edge sums by averaged quantities
involving the typical Delaunay edge. The local-time limit is obtained through
occupation-density arguments and Fourier representations of the local time,
in the spirit of \cite{Jaramillo21}. The treatment of the transition terms is
also inspired by the local-time approach developed in
\cite{Podolskij&Rosenbaum18} for fractional Brownian motion and by the power
variation analysis of Brown--Resnick processes in \cite{Robert20}. The main
additional difficulty here is the combination of a two-dimensional fractional
Brownian field, random Poisson observation sites, and the geometry of the
Delaunay graph.

The paper is organized as follows. Section \ref{sec:preliminaries} recalls
the required facts on the local time of the difference of two isotropic
fractional Brownian fields and on the Poisson--Delaunay graph. Section
\ref{sec:mainresults} introduces the normalized squared increment sum along
Delaunay edges, gives its decomposition into a Gaussian-type part and a
transition part, and states the main convergence theorem. Section
\ref{Proofs_main_results} proves the two main ingredients: the convergence of
the transition contribution to the local-time limit and the negligibility of
the Gaussian-type contribution at the relevant scale.  Section~\ref{sec:conclusion} concludes by summarizing the asymptotic mechanism and the technical tools used in the proof.
Some technical results are deferred in the Supplementary material. Section \ref{sec:L2_local_time} provides a
self-contained \(L^2\)-Fourier representation of the local time, while
Section \ref{sec:variance_truncated_quadratic_variation} proves the
auxiliary variance bound for the truncated quadratic variation which is used
in the proof of Proposition \ref{prop:Brownian_parts}. Section
\ref{sec:intermediary_results} collects some technical intermediary lemmas on increment
correlations, Delaunay-edge probabilities and Gaussian comparison bounds used
in the proofs of Propositions \ref{prop:twotrajectories} and \ref{prop:Brownian_parts}.

\section{Preliminaries}
\label{sec:preliminaries}

In this section, we recall some standard facts and introduce the notation used
throughout the paper.

\subsection{Local time of the difference of two isotropic fractional Brownian fields}
\label{subsec:local_time}

An isotropic fractional Brownian field is a centered Gaussian random field
\((W(x))_{x\in\mathbb R^2}\) such that \(W(0)=0\) a.s. and
\begin{equation}
\operatorname{Cov}(W(x),W(y))
=
\frac{\sigma^2}{2}
\left(
\|x\|^{2H}
+
\|y\|^{2H}
-
\|y-x\|^{2H}
\right),
\label{eq:defcovariance}
\end{equation}
for some \(H\in(0,1)\) and \(\sigma^2>0\). Here \(\|\cdot\|\) denotes the
Euclidean norm on \(\mathbb R^2\). The parameter \(\sigma\) is the scale
parameter, while \(H\) is the Hurst parameter. The field is \(H\)-self-similar and has stationary increments. Here \(H\)-self-similarity means that, for every \(\lambda>0\), the processes \(\left(W(\lambda x)\right)_{x\in\mathbb R^2}\) and \(\left(\lambda^H W(x)\right)_{x\in\mathbb R^2}\) have the same finite-dimensional distributions, while stationary increments means that, for every \(x_0\in\mathbb R^2\), the law of \(\left(W(x+x_0)-W(x_0)\right)_{x\in\mathbb R^2}\) does not depend on \(x_0\). 

Let \(W^{(1)}\) and \(W^{(2)}\) be two independent and identically distributed
isotropic fractional Brownian fields on \(\mathbb R^2\), with covariance
function given by \eqref{eq:defcovariance}. Define their difference field by $W^{(2\backslash 1)}(x)=        W^{(2)}(x)-W^{(1)}(x)$, $x\in\mathbb R^2$.
The local time of \(W^{(2\backslash 1)}\) provides the appropriate occupation
density for measuring the amount of time, in the spatial sense, that the
difference field spends near a given level. Let
 \(\nu^{(2\backslash 1)}\) be the occupation measure of
\(W^{(2\backslash 1)}\) over \(\mathbf C=(-1/2,1/2]^2\), namely
\[
        \nu^{(2\backslash 1)}(A)
        =
        \int_{\mathbf C}
        \mathbb I\left[
        W^{(2\backslash 1)}(x)\in A
        \right]\,\mathrm{d}x,
\]
for every Borel set \(A\subset\mathbb R\). For \(s,t\in\mathbf C\),
\[
\Delta(s,t)
:=
\mathbb E\left[
\left(
W^{(2\backslash 1)}(s)-W^{(2\backslash 1)}(t)
\right)^2
\right]
=
2\sigma^2\|s-t\|^{2H} .
\]
Since
\[
        \int_{\mathbf C}\Delta(s,t)^{-1/2}\,ds<\infty,
        \qquad t\in\mathbf C,
\]
it follows from the occupation density results of
\cite[Section 22]{Geman&Horowitz80} that the occupation measure
\(\nu^{(2\backslash 1)}\) admits a Lebesgue density. The local time of
\(W^{(2\backslash 1)}\) at level \(\ell\in\mathbb R\) is denoted by
\[
        L_{W^{(2\backslash 1)}}(\ell)
        =
        \frac{d\nu^{(2\backslash 1)}}{d\ell}(\ell).
\]
Equivalently, it is characterized by the occupation time formula
\[
        \int_{\mathbf C}
        g\left(W^{(2\backslash 1)}(x)\right)\,\mathrm{d}x
        =
        \int_{\mathbb R}
        g(\ell)L_{W^{(2\backslash 1)}}(\ell)\,d\ell,
\]
for every Borel measurable function \(g:\mathbb R\to\mathbb R\) for which the
two sides are well defined.

We shall also use two standard approximations of the local time. As recalled
in Section~\ref{sec:L2_local_time} of the Supplementary Material, for every \(\ell\in\mathbb R\),
\[
        L_{W^{(2\backslash 1)}}(\ell)
        =
        \lim_{\varepsilon\downarrow0}
        \int_{\mathbf C}
        \frac{1}{\sqrt{2\pi\varepsilon}}
        \exp\left(
        -
        \frac{
        \left(W^{(2\backslash 1)}(x)-\ell\right)^2
        }{2\varepsilon}
        \right)\,\mathrm{d}x
\]
in \(L^2\). Moreover, the Fourier representation
\begin{equation}
        L_{W^{(2\backslash 1)}}(\ell)
        =
        \frac{1}{2\pi}
        \lim_{M\to\infty}
        \int_{-M}^{M}
        \int_{\mathbf C}
        \exp\left(
        i\xi\left(W^{(2\backslash 1)}(x)-\ell\right)
        \right)\,\mathrm{d}x\,\mathrm{d}\xi
\label{eq:localtimeintegral}
\end{equation}
also holds in \(L^2\).

\subsection{Poisson--Delaunay graph}
\label{subsec:poisson_delaunay_graph}

Let \(P_N\) be a homogeneous Poisson point process with intensity \(N\) on
\(\mathbb R^2\). The Delaunay triangulation \(\operatorname{Del}(P_N)\) is the
triangulation with vertex set \(P_N\) such that the circumdisk of each triangle
contains no point of \(P_N\) in its interior. It is unique almost surely for a
Poisson point process; see, for instance, \cite[p.~478]{Schneider&Weil08}.
In the sequel, we only use the graph formed by the vertices and edges of this
triangulation, which we call the Poisson--Delaunay graph.

We recall the notion of typical cell for the Poisson--Delaunay tessellation
associated with a unit-intensity Poisson point process \(P_1\). With each cell
\(C\in\operatorname{Del}(P_1)\), associate its circumcenter \(z(C)\). If
\(\mathbf B\subset\mathbb R^2\) is a Borel set with area
\(a(\mathbf B)\in(0,\infty)\), the cell intensity is
\[
        \beta_2
        =
        \frac{1}{a(\mathbf B)}
        \mathbb E\left[
        \left|
        \{C\in\operatorname{Del}(P_1):z(C)\in\mathbf B\}
        \right|
        \right].
\]
It is well known that \(\beta_2=2\); see Theorem 10.2.9 in
\cite{Schneider&Weil08}. The typical cell \(\mathcal C\) is defined through
the identity
\[
        \mathbb E[g(\mathcal C)]
        =
        \frac{1}{\beta_2 a(\mathbf B)}
        \mathbb E\left[
        \sum_{C\in\operatorname{Del}(P_1):z(C)\in\mathbf B}
        g(C)
        \right],
\]
for every non-negative measurable translation-invariant function
\(g:\mathcal K_2\to\mathbb R\), where \(\mathcal K_2\) denotes the space of
compact convex subsets of \(\mathbb R^2\), endowed with the Fell topology (see \cite{Schneider&Weil08}). The
right-hand side does not depend on the particular choice of \(\mathbf B\), by
stationarity of \(P_1\) and translation invariance of \(g\).

The distribution of \(\mathcal C\) admits the following integral
representation; see Theorem 10.4.4 in \cite{Schneider&Weil08}:
\begin{equation}
        \mathbb E[g(\mathcal C)]
        =
        \frac{1}{6}
        \int_0^\infty
        \int_{(\mathbf S^1)^3}
        r^3 e^{-\pi r^2}
        a(\Delta(u_1,u_2,u_3))                        
        g(\Delta(ru_1,ru_2,ru_3))
        \,\sigma(du_1)\sigma(du_2)\sigma(du_3)\,dr .
\label{eq:typicalcell}
\end{equation}
Here \(\mathbf S^1\) is the unit circle of \(\mathbb R^2\),
\(\Delta(x_1,x_2,x_3)\) is the convex hull of \(\{x_1,x_2,x_3\}\), and
\(\sigma\) is the spherical Lebesgue measure on \(\mathbf S^1\), normalized by
\(\sigma(\mathbf S^1)=2\pi\). Equivalently, \(\mathcal C\) is equal in
distribution to
\[
        R\Delta(U_1,U_2,U_3),
\]
where \(R\) and \((U_1,U_2,U_3)\) are independent with densities
\[
        r\mapsto 2\pi^2 r^3 e^{-\pi r^2},
\]
and
\[
        (u_1,u_2,u_3)
        \mapsto
        \frac{a(\Delta(u_1,u_2,u_3))}{12\pi^2},
        \qquad (u_1,u_2,u_3)\in(\mathbf S^1)^3.
\]

We shall also use the distribution of the length of the typical
Poisson--Delaunay edge. The edge intensity of \(\operatorname{Del}(P_1)\) is $\beta_1=3$, again by Theorem 10.2.9 in \cite{Schneider&Weil08}. If \(D\) denotes the
length of the typical edge, then \(D\) has the same distribution as $R\|U_1-U_2\|$,
where \(R,U_1,U_2\) are obtained from the preceding typical-cell
representation. Its distribution function can be written as
\begin{align}
        \mathbb P(D\le \ell)
        &=
        \int_0^\ell f_D(d)\,dd \notag \\
        &=
        \frac{\pi}{3}
        \int_0^\infty
        \int_{(\mathbf S^1)^2}
        r^3e^{-\pi r^2}
        a(\Delta(u_1,u_2,e_1))
        \mathbb I\left[
        r\|u_1-u_2\|\le \ell
        \right]
        \,\sigma(du_1)\sigma(du_2)\,dr,
\label{eq:typicallength}
\end{align}
where \(e_1=(1,0)\), \(\ell>0\), and \(f_D\) denotes the density of \(D\).

Throughout the paper, we work with the one-skeleton of the Delaunay
triangulation. When two points \(x_1,x_2\in P_N\) are connected by an edge of
\(\operatorname{Del}(P_N)\), we say that they are Delaunay neighbors and write
\[
        x_1\sim x_2.
\]
Let \(\preceq\) denote the lexicographic order on \(\mathbb R^2\). For every
Borel set \(\mathbf B\subset\mathbb R^2\), define the oriented edge set
\[
        E_{N,\mathbf B}
        =
        \left\{
        (x_1,x_2)\in(P_N)^2:
        x_1\sim x_2,\,
        x_1\in\mathbf B,\,
        x_1\preceq x_2
        \right\}.
\]
We write
\[
        E_N:=E_{N,\mathbf C},
        \qquad
        \mathbf C=(-1/2,1/2]^2.
\]

We shall repeatedly use the following law of large numbers for the number of
Delaunay edges in the observation window:
\[
        \frac{|E_N|}{N}
        \xrightarrow{\mathrm{a.s.}}
        3,
        \qquad N\to\infty.
\]
This constant is the edge intensity of the stationary Poisson--Delaunay graph.
Such laws of large numbers follow from the general theory of stabilizing
geometric functionals of Poisson point processes; see, for example,
Chapter 10 in \cite{Schneider&Weil08}. 

\section{Main results}
\label{sec:mainresults}

Throughout this section, we assume that $H\in(0,1/2)$.
\(W^{(1)}\) and \(W^{(2)}\) are two independent and identically distributed
isotropic fractional Brownian fields on \(\mathbb R^2\), with covariance
function given by \eqref{eq:defcovariance} and their pointwise maximum is denoted by $W_\vee(x)=W^{(1)}(x)\vee W^{(2)}(x)$, $x\in\mathbb R^2$, while their difference is denoted by $W^{(2\backslash 1)}(x)=W^{(2)}(x)-W^{(1)}(x)$, $ x\in\mathbb R^2$.

For two distinct points \(x_1,x_2\in\mathbb R^2\), set
\[
        d_{1,2}=\|x_2-x_1\|.
\]
We introduce the normalized increment of the maximum field
\[
        U^{(W_\vee)}_{x_1,x_2}
        =
        \sigma^{-1}d_{1,2}^{-H}
        \left(
        W_\vee(x_2)-W_\vee(x_1)
        \right),
\]
for two distinct points $x_1$ and $x_2$.
The normalization \(\sigma d_{1,2}^{H}\) is the standard deviation of
the increment
\[
        W^{(i)}(x_2)-W^{(i)}(x_1),
        \qquad i=1,2,
\]
but it is not the standard deviation of the increment of \(W_\vee\). This
choice is nevertheless natural in the present infill regime, since away from
the random interface
\[
        \{x\in\mathbf C:W^{(1)}(x)=W^{(2)}(x)\},
\]
the field \(W_\vee\) locally coincides with one of the two underlying
fractional Brownian fields.

We consider the centered squared increment sum along the oriented Delaunay
edges in \(E_N\):
\[
        V^{(W_\vee)}_{2,N}
        =
        \frac{1}{\sqrt{|E_N|}}
        \sum_{(x_1,x_2)\in E_N}
        \left\{
        \left(U^{(W_\vee)}_{x_1,x_2}\right)^2-1
        \right\}.
\]
The centering constant \(1\) is the natural Gaussian centering, since for
\[
        U^{(i)}_{x_1,x_2}
        =
        \sigma^{-1}d_{1,2}^{-H}
        \left(
        W^{(i)}(x_2)-W^{(i)}(x_1)
        \right),
        \qquad i=1,2,
\]
one has
\[
        \mathbb E\left[
        \left(U^{(i)}_{x_1,x_2}\right)^2-1
        \right]
        =
        0.
\]
Strictly speaking, the statistics above are not defined on the events
\(\{|E_N|=0\}\). We shall use the convention that
\(V^{(W_\vee)}_{2,N}=0\) on \(\{|E_N|=0\}\). This convention is asymptotically immaterial, since these
exceptional events have exponentially small probability as \(N\to\infty\).

The first step is to decompose \(V^{(W_\vee)}_{2,N}\) according to whether the
same fractional Brownian field realizes the maximum at the initial point
\(x_1\), or whether a transition between the two fields may occur along the
edge \((x_1,x_2)\).

For any measurable function \(f:\mathbb R\to\mathbb R\), define
\begin{equation}
        \Psi_f(x,y,w)
        =
        \left(f(y+w)-f(x)\right)
        \mathbf 1_{\{x-y\le w\le 0\}}      
        +
        \left(f(x-w)-f(y)\right)
        \mathbf 1_{\{0\le w\le x-y\}} .
\end{equation}
Then, for every oriented edge \((x_1,x_2)\in E_N\),
\begin{equation}
\begin{aligned}
        f\left(U^{(W_\vee)}_{x_1,x_2}\right)
        &=
        f\left(U^{(1)}_{x_1,x_2}\right)
        \mathbf 1_{\{W^{(2\backslash 1)}(x_1)<0\}}       
        +
        f\left(U^{(2)}_{x_1,x_2}\right)
        \mathbf 1_{\{W^{(2\backslash 1)}(x_1)>0\}}       \\
        &\quad+
        \Psi_f\left(
        U^{(1)}_{x_1,x_2},
        U^{(2)}_{x_1,x_2},
        \frac{W^{(2\backslash 1)}(x_1)}
             {\sigma d_{1,2}^{H}}
        \right).
\end{aligned}
\label{eq:decompos_fU}
\end{equation}
Indeed, if \(W^{(2\backslash 1)}(x_1)<0\), then the maximum is realized by
\(W^{(1)}\) at \(x_1\), and a transition may occur at \(x_2\) only if
\[
        U^{(1)}_{x_1,x_2}-U^{(2)}_{x_1,x_2}
        \le
        \frac{W^{(2\backslash 1)}(x_1)}
             {\sigma d_{1,2}^{H}}
        \le 0.
\]
The case \(W^{(2\backslash 1)}(x_1)>0\) is symmetric.

We now apply \eqref{eq:decompos_fU} to the second Hermite polynomial $P^H_2(u)=u^2-1$, $u\in\mathbb R$, and write
$\Psi=\Psi_{P^H_2}$.
This gives the decomposition
\begin{equation}
        V^{(W_\vee)}_{2,N}
        =
        V^{(1)}_{2,N}
        +
        V^{(2)}_{2,N}
        +
        V^{(2/1)}_{2,N},
\label{eq:decompositionV}
\end{equation}
where
\[
        V^{(1)}_{2,N}
        =
        \frac{1}{\sqrt{|E_N|}}
        \sum_{(x_1,x_2)\in E_N}
        \left\{
        \left(U^{(1)}_{x_1,x_2}\right)^2-1
        \right\}
        \mathbf 1_{\{W^{(2\backslash 1)}(x_1)<0\}},
\]
\[
        V^{(2)}_{2,N}
        =
        \frac{1}{\sqrt{|E_N|}}
        \sum_{(x_1,x_2)\in E_N}
        \left\{
        \left(U^{(2)}_{x_1,x_2}\right)^2-1
        \right\}
        \mathbf 1_{\{W^{(2\backslash 1)}(x_1)>0\}},
\]
and
\[
        V^{(2/1)}_{2,N}
        =
        \frac{1}{\sqrt{|E_N|}}
        \sum_{(x_1,x_2)\in E_N}
        \Psi\left(
        U^{(1)}_{x_1,x_2},
        U^{(2)}_{x_1,x_2},
        \frac{W^{(2\backslash 1)}(x_1)}
             {\sigma d_{1,2}^{H}}
        \right).
\]
The first two terms correspond to edges for which the same field realizes the
maximum at the initial point \(x_1\). The last term is the transition term. It
is nonzero only when \(W^{(2\backslash 1)}(x_1)\) is of the same order as the
typical increment scale \(d_{1,2}^{H}\).

We now introduce the deterministic function which appears in the local-time
limit. Let us recall that \(D\) denotes the length of the typical Poisson--Delaunay edge, with
density \(f_D\) defined in \eqref{eq:typicallength}. For \(z\in\mathbb R\),
set
\[
        F_2(z)
        =
        \int_{\mathbb R^2\times\mathbb R_+}
        \Psi
        \left(
        x,y,\frac{z}{d^{H}}
        \right)
        \frac{1}{2\pi}
        e^{-(x^2+y^2)/2}
        f_D(d)\,\mathrm{d}x\,\mathrm{d}y\,\mathrm{d}d.
\]
The constant in the limiting theorem is
\[
        c_{V_2}
        =
        \int_{\mathbb R}F_2(z)\,\mathrm{d}z.
\]
The integrability of \(F_2\), and in fact its rapid decay, is verified in the
proof of Proposition~\ref{prop:twotrajectories}.

\begin{proposition}[Limit of the transition contribution]
\label{prop:twotrajectories}
Let \(W^{(1)}\) and \(W^{(2)}\) be two independent and identically distributed
isotropic fractional Brownian fields with covariance function
\eqref{eq:defcovariance}. Assume that \(H\in(0,1/2)\). Then, as
\(N\to\infty\),
\[
        \frac{\sqrt 3}{3}
        N^{-(2-2H)/4}
        V^{(2/1)}_{2,N}
        \xrightarrow{L^2}
        c_{V_2}
        L_{W^{(2\backslash 1)}}(0).
\]
\end{proposition}

The factor \(\sqrt 3/3\) comes from the normalization by
\(\sqrt{|E_N|}\) and from the edge-intensity law of large numbers: $|E_N|/N\xrightarrow{\mathrm{a.s.}}3$.

The second ingredient shows that the two Gaussian-type terms in
\eqref{eq:decompositionV} are negligible at the scale of the transition
contribution.

\begin{proposition}[Negligibility of the Gaussian-type contribution]
\label{prop:Brownian_parts}
Let \(W^{(1)}\) and \(W^{(2)}\) be two independent and identically distributed
isotropic fractional Brownian fields with covariance function
\eqref{eq:defcovariance}. Assume that \(H\in(0,1/2)\). Then, as
\(N\to\infty\),
\[
        N^{-(2-2H)/4}
        \left(
        V^{(1)}_{2,N}+V^{(2)}_{2,N}
        \right)
        \xrightarrow{L^2}
        0.
\]
\end{proposition}

Combining Propositions~\ref{prop:twotrajectories} and
\ref{prop:Brownian_parts} yields the main theorem.

\begin{theorem}[Limit theorem for the pointwise maximum]
\label{thm:main}
Under the assumptions of Propositions~\ref{prop:twotrajectories} and
\ref{prop:Brownian_parts}, as \(N\to\infty\),
\[
        \frac{\sqrt 3}{3}
        N^{-(2-2H)/4}
        V^{(W_\vee)}_{2,N}
        \xrightarrow{L^2}
        c_{V_2}
        L_{W^{(2\backslash 1)}}(0).
\]
\end{theorem}

Theorem~\ref{thm:main} shows that the asymptotic behavior of
\(V^{(W_\vee)}_{2,N}\) is governed by the transition contribution
\(V^{(2/1)}_{2,N}\). More precisely, the function \(\Psi\) is nonzero only
when
\[
        \left|
        W^{(2\backslash 1)}(x_1)
        \right|
        \lesssim
        \sigma d_{1,2}^{H}
        \left|
        U^{(1)}_{x_1,x_2}-U^{(2)}_{x_1,x_2}
        \right|,
\]
where the notation \(\lesssim\) is used in the usual sense: the left-hand side is bounded above by the right-hand side up to a multiplicative constant independent of \(N\) and of the edge.
Since a typical Poisson--Delaunay edge in the fixed window has length of order
\(N^{-1/2}\), the relevant transition band has width of order
\(N^{-H/2}\) around the random interface
$\{x\in\mathbf C:W^{(2\backslash 1)}(x)=0\}$.
The occupation density formula suggests the heuristic approximation
\[
        \int_{\mathbf C}
        \mathbf 1_{\{|W^{(2\backslash 1)}(x)|\le N^{-H/2}\}}\,\mathrm{d}x
        \approx
        2N^{-H/2}
        L_{W^{(2\backslash 1)}}(0).
\]
Thus the number of Delaunay edges whose initial point lies in the transition
band is of order $N\times N^{-H/2}=N^{1-H/2}$.
After the normalization by \(\sqrt{|E_N|}\), which is of order \(N^{1/2}\),
the transition contribution has size $N^{1-H/2}N^{-1/2}=N^{(2-2H)/4}$.
This explains both the normalization in Theorem~\ref{thm:main} and the
appearance of the local time \(L_{W^{(2\backslash 1)}}(0)\) in the limit.
The convergence is in $L^2$ because the leading random fluctuations are
averaged out conditionally on the interface, leaving the occupation-density
functional as the limiting object.

\section{Proofs}
\label{Proofs_main_results}

\subsection{Proof of Proposition \protect\ref{prop:twotrajectories}}

Without loss of generality, we assume that \(\sigma=1\). We write
\[
        X(x)=W^{(2\backslash 1)}(x),
        \qquad x\in\mathbb R^2.
\]

The proof follows the idea that the transition contribution behaves, after
averaging over the local Delaunay geometry, like a discrete approximation of
the occupation density of the difference field \(X\) at zero. More precisely,
the transition term is first rewritten as an edge functional depending on the
rescaled value \(X(x_1)/d_{1,2}^{H}\). Averaging the local edge
configuration around \(x_1\) gives a deterministic kernel \(F_f\), evaluated at
the rescaled field \(N^{H/2}X(x_1)\). Since the typical Delaunay edge has
length of order \(N^{-1/2}\), the factor \(N^{H/2}\) corresponds exactly
to the inverse width of the transition band around the interface
\(\{X=0\}\). The proof therefore consists in showing that this averaged
quantity converges in \(L^2\) to a constant times the local time \(L_X(0)\),
and then in proving that the original edge sum can be replaced by its averaged
version with an \(L^2\)-negligible error. The first step is based on the
Fourier representation of the local time and the Slivnyak--Mecke formula,
whereas the second step relies on residual covariance estimates for pairs of
Delaunay edges.

\strut

For any measurable function \(f:\mathbb R\to\mathbb R\), define
\[
        G_N^{(2\backslash 1)}[f]
        =
        \frac{1}{3}N^{H/2-1}
        \sum_{(x_1,x_2)\in E_N}
        \Psi_f\left(
        U^{(1)}_{x_1,x_2},
        U^{(2)}_{x_1,x_2},
        \frac{X(x_1)}{d_{1,2}^{H}}
        \right),
\]
and
\[
        G_{N,*}^{(2\backslash 1)}[f]
        =
        N^{H/2-1}
        \sum_{x\in P_N\cap\mathbf C}
        F_f\left(N^{H/2}X(x)\right),
\]
where
\[
        F_f(z)
        =
        \int_{\mathbb R^2\times\mathbb R_+}
        \Psi_f\left(x,y,\frac{z}{d^{H}}\right)
        \frac{1}{2\pi}e^{-(x^2+y^2)/2}f_D(d)\,\mathrm{d}x\,\mathrm{d}y\,\mathrm{d}d.
\]

For \(f=P^H_2\), one has
\[
        G_N^{(2\backslash 1)}[P^H_2]
        =
        \frac{1}{3}N^{H/2-1}\sqrt{|E_N|}
        V^{(2/1)}_{2,N}.
\]
Equivalently,
\[
        \frac{\sqrt 3}{3}
        N^{-(2-2H)/4}
        V^{(2/1)}_{2,N}
        =
        \left(\frac{3N}{|E_N|}\right)^{1/2}
        G_N^{(2\backslash 1)}[P^H_2].
\]
By the law of large numbers for Poisson--Delaunay functionals,
\[
        \frac{|E_N|}{N}\to3
        \qquad\text{a.s. and in }L^p
\]
for every finite \(p\). In particular,
\[
        \left(\frac{3N}{|E_N|}\right)^{1/2}
        \to1
         \qquad\text{a.s.}
\]
and this
multiplicative factor is uniformly bounded in \(L^p\) for some \(p\ge2\).
Therefore, once we prove
\[
        G_N^{(2\backslash 1)}[P^H_2]
        \xrightarrow{L^2}
        c_{V_2}L_X(0),
\]
the same \(L^2\)-convergence holds for the normalized transition contribution
appearing in Proposition~\ref{prop:twotrajectories}.



We prove the result in three steps. First, we show that, for any \(f\) such
that \(F_f\) belongs to the Schwartz space $\mathcal S(\mathbb R)$,
\begin{equation}
        G_{N,*}^{(2\backslash 1)}[f]
        \xrightarrow{L^2}
        c_f L_X(0),
        \qquad
        c_f:=\int_{\mathbb R}F_f(z)\,\mathrm{d}z .
\label{eq:twotrajectoriesstep1}
\end{equation}
Second, we verify that \(F_{P^H_2}\) has the integrability and decay properties
needed to extend the first step from Schwartz functions to  \(f=P^H_2\).
Third, we prove that
\begin{equation}
        G_N^{(2\backslash 1)}[P^H_2]
        -
        G_{N,*}^{(2\backslash 1)}[P^H_2]
        \xrightarrow{L^2}0.
\label{eq:twotrajectoriesstep2}
\end{equation}

\paragraph{(i) Proof of the \(L^2\)-convergence in \eqref{eq:twotrajectoriesstep1}.}

Let \(f\) be such that \(F_f\in\mathcal S(\mathbb R)\). For simplicity, write
\[
        F=F_f,
        \qquad
        c=c_f=\int_{\mathbb R}F(z)\,\mathrm{d}z.
\]
We prove that
\[
        \mathbb E\left[
        \left(
        G_{N,*}^{(2\backslash 1)}[f]-cL_X(0)
        \right)^2
        \right]
        \longrightarrow 0.
\]
Equivalently, it is enough to identify the limits of the three terms
\[
        \mathbb E\left[
        \left(G_{N,*}^{(2\backslash 1)}[f]\right)^2
        \right],
        \qquad
        \mathbb E\left[
        G_{N,*}^{(2\backslash 1)}[f]L_X(0)
        \right],
        \qquad
        \mathbb E[L_X(0)^2].
\]

Let \(\Sigma_{x,x'}\) denote the covariance matrix of $(X(x),X(x'))$.
By the Fourier representation of the local time in
\eqref{eq:localtimeintegral},
\begin{equation}
        \mathbb E\left[
        \left(cL_X(0)\right)^2
        \right]
        =
        \frac{c^2}{(2\pi)^2}
        \int_{\mathbb R^2}
        \int_{\mathbf C^2}
        \exp\left(
        -\frac12
        \vec \xi^{\,\top}\Sigma_{x,x'}\vec \xi
        \right)
        \,\mathrm{d}x\,\mathrm{d}x'\,\mathrm{d}\vec \xi,
\label{eq:termL2}
\end{equation}
where \(\vec\xi=(\xi,\xi')^\top\) and \(\mathrm{d}\vec\xi=\mathrm{d}\xi\,\mathrm{d}\xi'\). The
integrability of the right-hand side follows from Section
\ref{app:L2_local_time}.

We now study the second moment of \(G_{N,*}^{(2\backslash 1)}[f]\). Define
\[
        H_N(w)
        =
        N^{H/2}F(N^{H/2}w),
        \qquad w\in\mathbb R.
\]
Then
\[
        G_{N,*}^{(2\backslash 1)}[f]
        =
        \frac{1}{N}
        \sum_{x\in P_N\cap\mathbf C}
        H_N(X(x)).
\]
Since \(F\in\mathcal S(\mathbb R)\), Fourier inversion gives
\begin{equation}
        H_N(w)
        =
        \frac{1}{2\pi}
        \int_{\mathbb R}
        \int_{\mathbb R}
        F(y)
        \exp\left(
        i\xi\left(w-\frac{y}{N^{H/2}}\right)
        \right)
        \,\mathrm{d}y\,\mathrm{d}\xi .
\label{eq:HN_fourier}
\end{equation}

By the Slivnyak--Mecke formula,

\begin{equation}
        \mathbb E\left[
        \left(G_{N,*}^{(2\backslash 1)}[f]\right)^2
        \right]
        =
        \int_{\mathbf C^2}
        \mathbb E\left[
        H_N(X(x))H_N(X(x'))
        \right]\,\mathrm{d}x\,\mathrm{d}x'      
        +
        \frac{1}{N}
        \int_{\mathbf C}
        \mathbb E\left[
        H_N(X(x))^2
        \right]\,\mathrm{d}x .
\label{eq:second_moment_split}
\end{equation}
The second term is negligible. First, note that 
\(
        \|H_N\|_{L^2(\mathbb R)}^2
        =
        N^{H/2}\|F\|_{L^2(\mathbb R)}^2.
\)
Moreover, for \(x\neq0\), the random variable \(X(x)\) has a Gaussian density
bounded by \(C\|x\|^{-H}\). Hence
\[
        \mathbb E[H_N(X(x))^2]
        \le
        C N^{H/2}\|x\|^{-H}.
\]
Since \(H<2\), the function \(x\mapsto \|x\|^{-H}\) is
integrable over \(\mathbf C\). Therefore
\begin{equation}
        \frac{1}{N}
        \int_{\mathbf C}
        \mathbb E[H_N(X(x))^2]\,\mathrm{d}x
        =
        O\left(N^{-1+H/2}\right)
        \longrightarrow 0.
\label{eq:diagonal_negligible}
\end{equation}

For the first term in \eqref{eq:second_moment_split}, we use
\eqref{eq:HN_fourier}. Since \(F\) is rapidly decreasing and since
\[
        \int_{\mathbb R^2}
        \int_{\mathbf C^2}
        \exp\left(
        -\frac12
        \vec \xi^{\,\top}\Sigma_{x,x'}\vec \xi
        \right)
        \,\mathrm{d}x\,\mathrm{d}x'\,d\vec \xi
        <\infty,
\]
Fubini's theorem and dominated convergence yield
\[
\begin{aligned}
        &\int_{\mathbf C^2}
        \mathbb E\left[
        H_N(X(x))H_N(X(x'))
        \right]\,\mathrm{d}x\,\mathrm{d}x'        \\
        &\quad=
        \frac{1}{(2\pi)^2}
        \int_{\mathbf C^2}
        \int_{\mathbb R^2}
        \int_{\mathbb R^2}
        F(y)F(y')
        \exp\left(
        -i\frac{\xi y+\xi'y'}{N^{H/2}}
        \right)             
        \exp\left(
        -\frac12
        \vec \xi^{\,\top}\Sigma_{x,x'}\vec \xi
        \right)
        \,\mathrm{d}y\,\mathrm{d}y'\,d\vec \xi\,\mathrm{d}x\,\mathrm{d}x'          \\
        &\quad\longrightarrow
        \frac{c^2}{(2\pi)^2}
        \int_{\mathbb R^2}
        \int_{\mathbf C^2}
        \exp\left(
        -\frac12
        \vec \xi^{\,\top}\Sigma_{x,x'}\vec \xi
        \right)
        \,\mathrm{d}x\,\mathrm{d}x'\,d\vec \xi .
\end{aligned}
\]
Together with \eqref{eq:diagonal_negligible}, this gives
\begin{equation}
        \mathbb E\left[
        \left(G_{N,*}^{(2\backslash 1)}[f]\right)^2
        \right]
        \longrightarrow
        \frac{c^2}{(2\pi)^2}
        \int_{\mathbb R^2}
        \int_{\mathbf C^2}
        \exp\left(
        -\frac12
        \vec \xi^{\,\top}\Sigma_{x,x'}\vec \xi
        \right)
        \,\mathrm{d}x\,\mathrm{d}x'\,d\vec \xi .
\label{eq:termG2}
\end{equation}

It remains to identify the mixed term. By the Slivnyak--Mecke formula,
\[
        \mathbb E\left[
        G_{N,*}^{(2\backslash 1)}[f]\,cL_X(0)
        \right]
        =
        c
        \int_{\mathbf C}
        \mathbb E\left[
        H_N(X(x))L_X(0)
        \right]\,\mathrm{d}x.
\]
Using \eqref{eq:HN_fourier} and the Fourier representation of \(L_X(0)\), we
obtain
\[
\begin{aligned}
        &\mathbb E\left[
        G_{N,*}^{(2\backslash 1)}[f]\,cL_X(0)
        \right]        \\
        &\quad=
        \frac{c}{(2\pi)^2}
        \int_{\mathbf C^2}
        \int_{\mathbb R^2}
        \int_{\mathbb R}
        F(y)
        \exp\left(
        -i\frac{\xi y}{N^{H/2}}
        \right)
        \exp\left(
        -\frac12
        \vec \xi^{\,\top}\Sigma_{x,x'}\vec \xi
        \right)
        \,\mathrm{d}y\,d\vec \xi\,\mathrm{d}x\,\mathrm{d}x' .
\end{aligned}
\]
The same domination as above gives
\begin{equation}
        \mathbb E\left[
        G_{N,*}^{(2\backslash 1)}[f]\,cL_X(0)
        \right]
        \longrightarrow
        \frac{c^2}{(2\pi)^2}
        \int_{\mathbb R^2}
        \int_{\mathbf C^2}
        \exp\left(
        -\frac12
        \vec \xi^{\,\top}\Sigma_{x,x'}\vec \xi
        \right)
        \,\mathrm{d}x\,\mathrm{d}x'\,d\vec \xi .
\label{eq:term_cross}
\end{equation}

Combining \eqref{eq:termL2}, \eqref{eq:termG2} and
\eqref{eq:term_cross}, we get
\[
        \mathbb E\left[
        \left(
        G_{N,*}^{(2\backslash 1)}[f]-cL_X(0)
        \right)^2
        \right]
        \longrightarrow 0.
\]
This proves \eqref{eq:twotrajectoriesstep1}.

\strut

\paragraph{(ii) Regularity and decay of \(F_{P^H_2}\)}

It remains to study the regularity and decay of \(F_{P^H_2}\). We first
compute explicitly the Gaussian average which enters its definition. Let
\(Z_1,Z_2\) be two independent standard Gaussian random variables and set
\[
        h(w)
        =
        \mathbb E\left[
        \Psi_{P^H_2}(Z_1,Z_2,w)
        \right],
        \qquad w\in\mathbb R.
\]
Then
\[
        F_{P^H_2}(z)
        =
        \int_0^\infty
        h\left(\frac{z}{d^{H}}\right)
        f_D(d)\,\mathrm{d}d .
\]
We now give an explicit expression for \(h\). By symmetry, \(h\) is even. For
\(w\ge0\), using the definition of \(\Psi_{P^H_2}\), we have
\[
        h(w)
        =
        \mathbb E\left[
        \left\{(Z_1-w)^2-Z_2^2\right\}
        \mathbf 1_{\{w\le Z_1-Z_2\}}
        \right].
\]
Introduce
\[
        A=Z_1-Z_2,
        \qquad
        B=Z_1+Z_2.
\]
Then \(A\) and \(B\) are independent centered Gaussian random variables with
variance \(2\), and
\[
        (Z_1-w)^2-Z_2^2
        =
        (A-w)(B-w).
\]
Hence, for \(w\ge0\),
\begin{equation*}
        h(w)
       =
        \mathbb E\left[
        (A-w)(B-w)\mathbf 1_{\{A\ge w\}}
        \right]       
        =
        -w\,
        \mathbb E\left[
        (A-w)\mathbf 1_{\{A\ge w\}}
        \right].
\end{equation*}
Since \(A\sim\mathcal N(0,2)\),
\[
        \mathbb E\left[
        A\mathbf 1_{\{A\ge w\}}
        \right]
        =
        \frac{1}{\sqrt{\pi}}e^{-w^2/4},
\]
and therefore, with \(a=|w|\),
\begin{equation}
        h(w)
        =
        -a\left[
        \frac{1}{\sqrt{\pi}}e^{-a^2/4}
        -
        a\,\overline\Phi\left(\frac{a}{\sqrt2}\right)
        \right],
        \qquad a=|w|.
\label{eq:h_explicit}
\end{equation}
In particular,
\[
        h(w)
        =
        -\frac{|w|}{\sqrt{\pi}}+O(w^2),
        \qquad w\to0.
\]
Consequently,
\begin{equation}
        F_{P^H_2}(z)
        =
        -\frac{|z|}{\sqrt{\pi}}
        \int_0^\infty d^{-H}f_D(d)\,\mathrm{d}d
        +
        O(z^2),
        \qquad z\to0.
\label{eq:F_H2_near_zero}
\end{equation}
The integral in \eqref{eq:F_H2_near_zero} is finite because the density of
the typical Delaunay edge satisfies \(f_D(d)=O(d^3)\) as \(d\downarrow0\).
Thus \(F_{P^H_2}\) is continuous and even, but it is not differentiable at
zero. In particular, \(F_{P^H_2}\) does not belong to the Schwartz space.

This lack of smoothness at the origin is harmless for the preceding argument.
What is needed is integrability and sufficiently fast decay at infinity. We
now verify this. From \eqref{eq:h_explicit}, for every \(k\ge0\) there exist
constants \(C_k,c_k>0\) such that, for \(w>0\),
\[
        |h^{(k)}(w)|
        \le
        C_k(1+w)^{m_k}e^{-c_k w^2},
\]
for some integer \(m_k\ge0\). Moreover, by \eqref{eq:typicallength}, the
density \(f_D\) satisfies
\[
        f_D(d)\le C d^3e^{-\pi d^2},
        \qquad d>0.
\]
Therefore, for \(z>0\),
\[
\begin{aligned}
        |F_{P^H_2}^{(k)}(z)|
        &\le
        C_k
        \int_0^\infty
        d^{-kH}
        \left(1+\frac{z}{d^{H}}\right)^{m_k}
        \exp\left(
        -c_k\frac{z^2}{d^{2H}}
        \right)
        d^3e^{-\pi d^2}\,\mathrm{d}d .
\end{aligned}
\]
A standard Laplace bound gives
\[
        \frac{z^2}{d^{2H}}+d^2
        \ge
        c z^{4/(2+2H)},
        \qquad z\ge1,\ d>0.
\]
It follows that, for every \(k\ge0\), there exist constants \(C_k,c_k'>0\)
such that
\begin{equation}
        |F_{P^H_2}^{(k)}(z)|
        \le
        C_k
        \exp\left(
        -c_k' z^{4/(2+2H)}
        \right),
        \qquad z\ge1.
\label{eq:F_H2_decay}
\end{equation}
By evenness, the same estimate holds as \(z\to-\infty\). In particular,
\(F_{P^H_2}\in L^1(\mathbb R)\cap L^2(\mathbb R)\), and
\[
        c_{V_2}
        =
        \int_{\mathbb R}F_{P^H_2}(z)\,\mathrm{d}z
\]
is well defined.

The first step of the proof, initially stated for Schwartz functions, can
therefore be applied to \(F_{P^H_2}\) by a standard approximation argument. Let
\((F_m)_{m\ge1}\subset\mathcal S(\mathbb R)\) be a sequence such that
\(F_m\to F_{P^H_2}\) in \(L^1(\mathbb R)\cap L^2(\mathbb R)\), with the same
type of domination as in \eqref{eq:F_H2_decay}. Applying
\eqref{eq:twotrajectoriesstep1} to \(F_m\), and then letting \(m\to\infty\),
gives the desired convergence for \(F_{P^H_2}\).

\paragraph{(iii) Proof of the \(L^2\)-convergence in 
\eqref{eq:twotrajectoriesstep2}.}

We now prove that
\[
        G_N^{(2\backslash 1)}[P^H_2]
        -
        G_{N,*}^{(2\backslash 1)}[P^H_2]
        \xrightarrow{L^2}0.
\]
Throughout this part of the proof we write
\[
        X(x)=W^{(2\backslash 1)}(x),
        \qquad
        F=F_{P^H_2},
        \qquad
        \Psi=\Psi_{P^H_2}.
\]
For an oriented edge \(e=(x_1,x_2)\in E_N\), set
\[
        d_e=\|x_2-x_1\|,
        \qquad
        \Psi_e
        =
        \Psi\left(
        U^{(1)}_{x_1,x_2},
        U^{(2)}_{x_1,x_2},
        \frac{X(x_1)}{d_e^{H}}
        \right),
\]
and
\[
        F_{N}(x_1)
        =
        F\left(N^{H/2}X(x_1)\right).
\]
Then
\[
        G_N^{(2\backslash 1)}[P^H_2]
        =
        \frac{1}{3}N^{H/2-1}
        \sum_{e=(x_1,x_2)\in E_N}\Psi_e,
\]
whereas
\[
        G_{N,*}^{(2\backslash 1)}[P^H_2]
        =
        N^{H/2-1}
        \sum_{x\in P_N\cap\mathbf C}F_N(x).
\]

We split the difference into two terms:
\[
\begin{aligned}
        G_N^{(2\backslash 1)}[P^H_2]
        -
        G_{N,*}^{(2\backslash 1)}[P^H_2]
        &=
        \frac{1}{3}R_N
        +
        A_N,
\end{aligned}
\]
where
\[
        R_N
        =
        N^{H/2-1}
        \sum_{e=(x_1,x_2)\in E_N}
        \left(\Psi_e-F_N(x_1)\right),
\]
and
\[
        A_N
        =
        N^{H/2-1}
        \left\{
        \frac{1}{3}
        \sum_{e=(x_1,x_2)\in E_N}F_N(x_1)
        -
        \sum_{x\in P_N\cap\mathbf C}F_N(x)
        \right\}.
\]
We shall prove that
\[
        A_N\xrightarrow{L^2}0,
        \qquad
        R_N\xrightarrow{L^2}0.
\]

\medskip
We first handle the replacement of the vertex sum by the edge sum.

\paragraph{Step 1: $A_N\xrightarrow{L^2}0$.}

For \(x\in P_N\cap\mathbf C\), let
\[
        D_N^+(x)
        =
        \#\{y\in P_N:\ (x,y)\in E_N\}
\]
be the number of oriented Delaunay edges in \(E_N\) starting from \(x\).
Then
\[
        \frac{1}{3}
        \sum_{e=(x_1,x_2)\in E_N}F_N(x_1)
        -
        \sum_{x\in P_N\cap\mathbf C}F_N(x)
        =
        \sum_{x\in P_N\cap\mathbf C}
        \left(
        \frac{D_N^+(x)}{3}-1
        \right)F_N(x).
\]
The score
\[
        \xi_N(x,P_N)
        =
        \frac{D_N^+(x)}{3}-1
\]
is a translation-invariant score for the
Poisson--Delaunay graph. Moreover, by the edge-intensity identity
\(\beta_1=3\),
\[
        \mathbb E^0[\xi_N(0,P_N)]=0,
\]
where \(\mathbb E^0\) denotes expectation under the Palm distribution.
Consequently, the standard second-moment estimate for Poisson functionals with centered scores gives
\begin{equation}
\label{eq:AN_second_moment}
        \mathbb E[A_N^2]
        \le
        C N^{H-2}
        \left\{
        N\int_{\mathbf C}
        \mathbb E\left[F_N(x)^2\right]\,\mathrm{d}x
        +
        N^2
        \int_{\mathbf C^2}
        e^{-c\sqrt N\|x-y\|}
        \left(
        \mathbb E[F_N(x)^2]\,
        \mathbb E[F_N(y)^2]
        \right)^{1/2}
        \mathrm{d}x\,\mathrm{d}y
        \right\}.
\end{equation}
The first term in braces corresponds to the diagonal contribution. The second
one controls the off-diagonal contribution.

The second term in \eqref{eq:AN_second_moment} is of the same
order as the first one. Indeed, by Young's inequality and the fact that
\[
        \int_{\mathbb R^2}e^{-c\sqrt N\|z\|}\,\mathrm{d}z
        =
        O(N^{-1}),
\]
we have
\begin{equation*}
        N^2
        \int_{\mathbf C^2}
        e^{-c\sqrt N\|x-y\|}
        \left(
        \mathbb E[F_N(x)^2]\,
        \mathbb E[F_N(y)^2]
        \right)^{1/2}
        \mathrm{d}x\,\mathrm{d}y                                     
        \le
        C N
        \int_{\mathbf C}
        \mathbb E[F_N(x)^2]\,\mathrm{d}x .
\end{equation*}
Thus
\begin{equation}
\label{eq:AN_second_moment_reduced}
        \mathbb E[A_N^2]
        \le
        C N^{H-1}
        \int_{\mathbf C}
        \mathbb E\left[F_N(x)^2\right]\,\mathrm{d}x .
\end{equation}

Since \(F\) is rapidly decreasing, for \(x\neq0\),
\[
        \mathbb E\left[
        F\left(N^{H/2}X(x)\right)^2
        \right]
        \le
        C N^{-H/2}\|x\|^{-H}.
\]
Since \(x\mapsto \|x\|^{-H}\) is integrable over \(\mathbf C\), we get
\[
        \int_{\mathbf C}
        \mathbb E\left[F_N(x)^2\right]\,\mathrm{d}x
        \le
        C N^{-H/2}.
\]
Combining this bound with \eqref{eq:AN_second_moment_reduced}, we obtain
\[
        \mathbb E[A_N^2]
        \le
        C N^{H-1}N^{-H/2}
        =
        C N^{-1+H/2}.
\]
Since \(2H<1\), this tends to \(0\). Hence
\[
        A_N\xrightarrow{L^2}0.
\]

It remains to prove that \(R_N\to0\) in $L^2$. 

\paragraph{Step 2: \(R_N\xrightarrow{L^2}0\).}

Recall that
\[
        R_N
        =
        N^{H/2-1}
        \sum_{e=(x_1,x_2)\in E_N}
        \left(\Psi_e-F_N(x_1)\right),
\]
where
\[
        F_N(x)=F\left(N^{H/2}X(x)\right)
\]
and
\[
        \Psi_e
        =
        \Psi\left(
        U^{(1)}_{x_1,x_2},
        U^{(2)}_{x_1,x_2},
        \frac{X(x_1)}{\|x_2-x_1\|^{H}}
        \right).
\]
For \(e=(x_1,x_2)\in E_N\), define
\[
        r_{N,e}
        =
        \Psi_e-F_N(x_1).
\]
Then
\begin{equation}
        \mathbb E[R_N^2]
        =
        N^{H-2}
        \mathbb E\left[
        \sum_{e,e'\in E_N}
        r_{N,e}r_{N,e'}
        \right].
\label{eq:RN_second_moment_start}
\end{equation}

We first remove the local pairs of edges. Let
\[
        \mathcal L_N
        =
        \{(e,e')\in E_N^2:\ e\cap e'\neq\varnothing\}
\]
be the set of ordered pairs of oriented edges sharing at least one endpoint.
We claim that
\begin{equation}
        N^{H-2}
        \mathbb E\left[
        \sum_{(e,e')\in\mathcal L_N}
        |r_{N,e}r_{N,e'}|
        \right]
        \longrightarrow0.
\label{eq:local_pairs_negligible_RN}
\end{equation}
Indeed, the elementary bound
\[
        |\Psi(u,v,w)|
        \le
        C(1+u^2+v^2)
        \mathbf 1_{\{|w|\le |u-v|\}}
\]
implies
\[
        \mathbb E[\Psi_e^2]\le C
\]
uniformly in the edge \(e\). Moreover, \(F\) is bounded, since it has
stretched-exponential decay. Hence
\[
        \mathbb E[r_{N,e}^2]\le C
\]
uniformly in \(N\) and \(e\). By Cauchy's inequality,
\[
        \mathbb E[|r_{N,e}r_{N,e'}|]\le C.
\]
The expected number of ordered adjacent pairs of Delaunay edges in
\(\mathbf C\) is \(O(N)\), because the degree of the typical Poisson--Delaunay
vertex has finite moments of all orders. Therefore
\[
        N^{H-2}
        \mathbb E\left[
        \sum_{(e,e')\in\mathcal L_N}
        |r_{N,e}r_{N,e'}|
        \right]
        \le
        C N^{H-2}N
        =
        C N^{-1+H}
        \longrightarrow0,
\]
since \(2H<1\). This proves
\eqref{eq:local_pairs_negligible_RN}.

It remains to treat the contribution of pairs of disjoint edges. Let
\[
        \mathcal D_N
        =
        \{(e,e')\in E_N^2:\ e\cap e'=\varnothing\}.
\]
For two disjoint oriented edges
\[
        e=(x_1,x_2),
        \qquad
        e'=(x_3,x_4),
\]
write
\[
        r_{N;x_1,x_2}
        =
        \Psi\left(
        U^{(1)}_{x_1,x_2},
        U^{(2)}_{x_1,x_2},
        \frac{X(x_1)}{\|x_2-x_1\|^{H}}
        \right)
        -
        F\left(N^{H/2}X(x_1)\right).
\]
Let \(p_{2,N}(x_1,x_2,x_3,x_4)\) denote the probability that
\(x_1\sim x_2\) and \(x_3\sim x_4\), with the prescribed orientation, when
the four fixed points are inserted into the Poisson point process. By the
Slivnyak--Mecke formula, the contribution of disjoint pairs in
\eqref{eq:RN_second_moment_start} is
\begin{equation}
        \mathcal R_N
        =
        N^{H+2}
        \int_{\mathbf C^2}
        \int_{(\mathbb R^2)^2}
        \mathbb E\left[
        r_{N;x_1,x_2}
        r_{N;x_3,x_4}
        \right]                                      
        p_{2,N}(x_1,x_2,x_3,x_4)
        \,\mathrm{d}x_2\,\mathrm{d}x_4\,\mathrm{d}x_1\,\mathrm{d}x_3 .
\label{eq:RN_disjoint_contribution}
\end{equation}
We shall prove that $\mathcal R_N\longrightarrow0$.

The key estimate is the following residual covariance kernel bound.

\begin{lemma}[Residual covariance kernel]
\label{lem:residual_covariance_kernel}
For \(x_1,x_3\in\mathbf C\), set
\[
        K_N(x_1,x_3)
        =
        N^{2+H}
        \int_{(\mathbb R^2)^2}
        \mathbb E
        \left[
        r_{N;x_1,x_2}
        r_{N;x_3,x_4}
        \right]
        p_{2,N}(x_1,x_2,x_3,x_4)
        \,\mathrm{d}x_2\,\mathrm{d}x_4 .
\]
Let
\[
        \Gamma_{x_1,x_3}
        =
        \operatorname{Cov}\bigl(X(x_1),X(x_3)\bigr).
\]
Then there exists \(C<\infty\), independent of \(N,x_1,x_3\), such that
\begin{equation}
        |K_N(x_1,x_3)|
        \le
        \frac{C}{\sqrt{\det\Gamma_{x_1,x_3}}},
        \qquad x_1\neq x_3.
\label{eq:kernel_bound_det}
\end{equation}
Consequently,
\begin{equation}
        \lim_{\delta\downarrow0}
        \sup_{N\ge1}
        \int_{\mathbf C^2}
        \mathbf 1_{\{\|x_1-x_3\|\le\delta\}}
        |K_N(x_1,x_3)|\,\mathrm{d}x_1\mathrm{d}x_3
        =
        0 .
\label{eq:kernel_bound_near_integral}
\end{equation}
Moreover, for every \(\delta>0\),
\begin{equation}
        \int_{\mathbf C^2}
        \mathbf 1_{\{\|x_1-x_3\|>\delta\}}
        K_N(x_1,x_3)\,\mathrm{d}x_1\mathrm{d}x_3
        \longrightarrow0.
\label{eq:kernel_bound_far_integral}
\end{equation}
\end{lemma}

Using the definition of \(K_N\), the disjoint contribution
\eqref{eq:RN_disjoint_contribution} can be written simply as
\begin{equation}
        \mathcal R_N
        =
        \int_{\mathbf C^2}K_N(x_1,x_3)\,\mathrm{d}x_1\,\mathrm{d}x_3.
\label{eq:RN_as_integral_KN}
\end{equation}
Let \(\delta>0\). We split
\[
        \mathcal R_N
        =
        \mathcal R_{N,\delta}^{(1)}
        +
        \mathcal R_{N,\delta}^{(2)},
\]
where
\[
        \mathcal R_{N,\delta}^{(1)}
        =
        \int_{\mathbf C^2}
        \mathbf 1_{\{\|x_1-x_3\|\le\delta\}}
        K_N(x_1,x_3)\,\mathrm{d}x_1\,\mathrm{d}x_3,
\]
and
\[
        \mathcal R_{N,\delta}^{(2)}
        =
        \int_{\mathbf C^2}
        \mathbf 1_{\{\|x_1-x_3\|>\delta\}}
        K_N(x_1,x_3)\,\mathrm{d}x_1\,\mathrm{d}x_3.
\]

By \eqref{eq:kernel_bound_near_integral}, we have
\begin{equation}
        \lim_{\delta\downarrow0}
        \sup_{N\ge1}
        |\mathcal R_{N,\delta}^{(1)}|
        =
        0.
\label{eq:RN_near_diagonal_vanishes}
\end{equation}
More explicitly, using the determinant bound in
\eqref{eq:kernel_bound_det} and the estimate from
Section~\ref{app:L2_local_time}, one obtains, for some \(C<\infty\),
\[
        \sup_{N\ge1}
        |\mathcal R_{N,\delta}^{(1)}|
        \le
        C\delta^{2-H}.
\]
On the other hand, by \eqref{eq:kernel_bound_far_integral},
\begin{equation}
        \mathcal R_{N,\delta}^{(2)}
        \longrightarrow0
        \qquad
        \text{for every fixed }\delta>0.
\label{eq:RN_far_diagonal_vanishes}
\end{equation}
Combining
\eqref{eq:RN_near_diagonal_vanishes} and
\eqref{eq:RN_far_diagonal_vanishes}, we get
\[
        \mathcal R_N\longrightarrow0.
\]

Together with the negligible contribution of local pairs in
\eqref{eq:local_pairs_negligible_RN}, this gives $\mathbb E[R_N^2]\longrightarrow0$.
Hence $R_N\xrightarrow{L^2}0$.

Combining this with the replacement estimate \(A_N\xrightarrow{L^2}0\), we
obtain
\[
        G_N^{(2\backslash 1)}[P^H_2]
        -
        G_{N,*}^{(2\backslash 1)}[P^H_2]
        \xrightarrow{L^2}0,
\]
which proves \eqref{eq:twotrajectoriesstep2}.

\strut

\begin{prooft}{Lemma \ref{lem:residual_covariance_kernel}}

We prove the three assertions in the statement. The proof is divided into two
parts. We first establish the uniform determinant bound
\eqref{eq:kernel_bound_det}, which immediately implies the near-diagonal
domination \eqref{eq:kernel_bound_near_integral}. We then prove the
far-field convergence \eqref{eq:kernel_bound_far_integral}.

\medskip

\noindent\textbf{Step 1: determinant bound and near-diagonal domination.}

Let
\[
        \Gamma_{x_1,x_3}
        =
        \operatorname{Cov}\bigl(X(x_1),X(x_3)\bigr)
\]
be the covariance matrix of the Gaussian vector \((X(x_1),X(x_3))\). We first
prove that
\[
        |K_N(x_1,x_3)|
        \le
        \frac{C}{\sqrt{\det\Gamma_{x_1,x_3}}},
        \qquad x_1\neq x_3.
\]

Recall that
\[
        r_{N;x_1,x_2}
        =
        \Psi\left(
        U^{(1)}_{x_1,x_2},
        U^{(2)}_{x_1,x_2},
        \frac{X(x_1)}{\|x_2-x_1\|^{H}}
        \right)
        -
        F\left(N^{H/2}X(x_1)\right).
\]
For notational simplicity, set
\[
        \ell=\|x_2-x_1\|,
        \qquad
        \ell'=\|x_4-x_3\|,
\]
and define
\[
        \Psi_e
        =
        \Psi\left(
        U^{(1)}_{x_1,x_2},
        U^{(2)}_{x_1,x_2},
        \frac{X(x_1)}{\ell^{H}}
        \right),
\]
\[
        \Psi_{e'}
        =
        \Psi\left(
        U^{(1)}_{x_3,x_4},
        U^{(2)}_{x_3,x_4},
        \frac{X(x_3)}{(\ell')^{H}}
        \right),
\]
and
\[
        F_N(x)=F\left(N^{H/2}X(x)\right).
\]
Then
\begin{equation}
        |r_{N;x_1,x_2}r_{N;x_3,x_4}|
        \le
        C\Big(
        |\Psi_e\Psi_{e'}|
        +
        |\Psi_e|\,|F_N(x_3)|             
        +
        |F_N(x_1)|\,|\Psi_{e'}|
        +
        |F_N(x_1)F_N(x_3)|
        \Big).
\label{eq:residual_product_split_kernel}
\end{equation}

We shall bound each term separately. We use the elementary estimate
\begin{equation}
        |\Psi(u,v,w)|
        \le
        C(1+u^2+v^2)\,
        \mathbf 1_{\{|w|\le |u-v|\}} .
\label{eq:Psi_basic_bound_kernel}
\end{equation}
Moreover, \(F\) has stretched-exponential decay: there exist constants
\(C,c,\beta>0\) and an integer \(q\ge0\) such that
\begin{equation}
        |F(z)|
        \le
        C(1+|z|)^q e^{-c|z|^\beta},
        \qquad z\in\mathbb R .
\label{eq:F_decay_kernel}
\end{equation}

We start with the term \(\mathbb E[|\Psi_e\Psi_{e'}|]\). Conditionally on
\((X(x_1),X(x_3))\), the vector
\[
        \left(
        U^{(1)}_{x_1,x_2},
        U^{(2)}_{x_1,x_2},
        U^{(1)}_{x_3,x_4},
        U^{(2)}_{x_3,x_4}
        \right)
\]
is Gaussian. Its conditional mean and covariance matrix are given by the
Gaussian regression formulas. 

By Lemma~\ref{lem:conditional_PsiPsi_bound}, conditionally on
$X(x_1)=z_1$, and $X(x_3)=z_3$,
we have
\begin{equation}
\begin{aligned}
        &\mathbb E\left[
        |\Psi_e\Psi_{e'}|
        \,\middle|\,
        X(x_1)=z_1,\ X(x_3)=z_3
        \right]                                           \\
        &\quad\le
        C\,
        \mathbb E\left[
        B_{e,e'}(z,\mathcal Z)^4
        \mathbf 1_{\{|z_1|\le C\ell^{H}B_{e,e'}(z,\mathcal Z)\}}
        \mathbf 1_{\{|z_3|\le C(\ell')^{H}B_{e,e'}(z,\mathcal Z)\}}
        \right],
\end{aligned}
\label{eq:PsiPsi_conditional_bound_kernel}
\end{equation}

where \(B_{e,e'}\) can be chosen as in
\eqref{eq:B_explicit_standardized}.

Integrating \eqref{eq:PsiPsi_conditional_bound_kernel} with respect to the
Gaussian density of \((X(x_1),X(x_3))\), we get
\begin{equation}
        \mathbb E\left[
        |\Psi_e\Psi_{e'}|
        \right]
        \le
        C
        \frac{
        \ell^{H}(\ell')^{H}
        }
        {\sqrt{\det\Gamma_{x_1,x_3}}}.
\label{eq:PsiPsi_unconditional_bound_kernel}
\end{equation}
Indeed, the two indicators in
\eqref{eq:PsiPsi_conditional_bound_kernel} restrict \((z_1,z_3)\) to a
rectangle with side lengths of order \(\ell^{H}\) and
\((\ell')^{H}\), up to polynomial Gaussian factors. Since the bivariate Gaussian density function satisfies
\[
        \phi_{\Gamma_{x_1,x_3}}(z_1,z_3)
        \le
        \frac{C}{\sqrt{\det\Gamma_{x_1,x_3}}},
\]
the integration over this rectangle yields the factor
\[
        \frac{
        \ell^{H}(\ell')^{H}
        }
        {\sqrt{\det\Gamma_{x_1,x_3}}},
\]
and the remaining polynomial factors have finite Gaussian moments.

The same argument, using the decay of \(F\), gives the mixed estimates
\begin{equation}
        \mathbb E\left[
        |\Psi_e|\,|F_N(x_3)|
        \right]
        \le
        C
        \frac{
        \ell^{H}N^{-H/2}
        }
        {\sqrt{\det\Gamma_{x_1,x_3}}},
\label{eq:PsiF_bound_kernel}
\end{equation}
and
\begin{equation}
        \mathbb E\left[
        |F_N(x_1)|\,|\Psi_{e'}|
        \right]
        \le
        C
        \frac{
        N^{-H/2}(\ell')^{H}
        }
        {\sqrt{\det\Gamma_{x_1,x_3}}}.
\label{eq:FPsi_bound_kernel}
\end{equation}
Finally, by the change of variables
\[
        u=N^{H/2}z_1,
        \qquad
        v=N^{H/2}z_3,
\]
and by the integrability of \(F\), we obtain
\begin{equation}
        \mathbb E\left[
        |F_N(x_1)F_N(x_3)|
        \right]
        \le
        C
        \frac{N^{-H}}
        {\sqrt{\det\Gamma_{x_1,x_3}}}.
\label{eq:FF_bound_kernel}
\end{equation}

Combining
\eqref{eq:residual_product_split_kernel},
\eqref{eq:PsiPsi_unconditional_bound_kernel},
\eqref{eq:PsiF_bound_kernel},
\eqref{eq:FPsi_bound_kernel} and
\eqref{eq:FF_bound_kernel}, we obtain
\begin{equation}
\begin{aligned}
        &\mathbb E\left[
        |r_{N;x_1,x_2}r_{N;x_3,x_4}|
        \right]                                                   \\
        &\quad\le
        \frac{C}{\sqrt{\det\Gamma_{x_1,x_3}}}
        \Big[
        \ell^{H}(\ell')^{H}
        +
        \ell^{H}N^{-H/2}
        +
        N^{-H/2}(\ell')^{H}
        +
        N^{-H}
        \Big].
\end{aligned}
\label{eq:residual_pair_length_bound_kernel}
\end{equation}

We now integrate with respect to the two-edge Delaunay probability. Using Lemma \ref{Le:estimatepN}, for every \(s,t\ge0\),
\begin{equation}
        \int_{(\mathbb R^2)^2}
        \ell^s(\ell')^t
        p_{2,N}(x_1,x_2,x_3,x_4)\,\mathrm{d}x_2\mathrm{d}x_4
        \le
        C_{s,t}N^{-2-(s+t)/2}.
\label{eq:Delaunay_pair_length_moment_kernel}
\end{equation}
Applying \eqref{eq:Delaunay_pair_length_moment_kernel} with 
\[
        (s,t)=\left(H,H\right),
        \qquad
        (s,t)=\left(H,0\right),
        \qquad
        (s,t)=\left(0,H\right),
        \qquad
        (s,t)=(0,0),
\]
we deduce from \eqref{eq:residual_pair_length_bound_kernel} that
\begin{equation*}
        \int_{(\mathbb R^2)^2}
        \mathbb E\left[
        |r_{N;x_1,x_2}r_{N;x_3,x_4}|
        \right]
        p_{2,N}(x_1,x_2,x_3,x_4)\,\mathrm{d}x_2\mathrm{d}x_4      
        \qquad\le
        C N^{-2-H}
        \frac{1}{\sqrt{\det\Gamma_{x_1,x_3}}}.
\end{equation*}
Multiplying by \(N^{2+H}\) yields
\[
        |K_N(x_1,x_3)|
        \le
        \frac{C}{\sqrt{\det\Gamma_{x_1,x_3}}},
\]
which proves \eqref{eq:kernel_bound_det}.

We now derive the near-diagonal domination. By the determinant estimate proved
in Section~\ref{app:L2_local_time},
\begin{equation}
        \frac{1}{\sqrt{\det\Gamma_{x_1,x_3}}}
        \le
        C
        \bigl(
        \|x_1\|^{-H}
        +
        \|x_3\|^{-H}
        \bigr)
        \|x_1-x_3\|^{-H}.
  \label{eq:bound_frac_1_sqrt(det)}      
\end{equation}
Hence
\[
\begin{aligned}
        &\int_{\mathbf C^2}
        \mathbf 1_{\{\|x_1-x_3\|\le\delta\}}
        |K_N(x_1,x_3)|\,\mathrm{d}x_1\mathrm{d}x_3                          \\
        &\qquad\le
        C
        \int_{\mathbf C^2}
        \mathbf 1_{\{\|x_1-x_3\|\le\delta\}}
        \bigl(
        \|x_1\|^{-H}
        +
        \|x_3\|^{-H}
        \bigr)
        \|x_1-x_3\|^{-H}
        \mathrm{d}x_1\mathrm{d}x_3 .
\end{aligned}
\]
The right-hand side is bounded by \(C\delta^{2-H}\), because both
\(\|x\|^{-H}\) and \(\|x_1-x_3\|^{-H}\) are locally integrable
in dimension two. Therefore
\[
        \lim_{\delta\downarrow0}
        \sup_{N\ge1}
        \int_{\mathbf C^2}
        \mathbf 1_{\{\|x_1-x_3\|\le\delta\}}
        |K_N(x_1,x_3)|\,\mathrm{d}x_1\mathrm{d}x_3
        =
        0,
\]
which proves \eqref{eq:kernel_bound_near_integral}.

\medskip

\noindent\textbf{Step 2: far-field convergence.}

We now prove \eqref{eq:kernel_bound_far_integral}. Fix \(\delta>0\). We first
work away from the origin. For \(\varepsilon>0\), set
\[
        \mathbf C_\varepsilon
        =
        \{x\in\mathbf C:\|x\|\ge\varepsilon\}.
\]
We shall prove that
\begin{equation}
        \sup_{\substack{x_1,x_3\in\mathbf C_\varepsilon\\
        \|x_1-x_3\|>\delta}}
        |K_N(x_1,x_3)|
        \longrightarrow0.
\label{eq:kernel_far_eps_proof}
\end{equation}
Once \eqref{eq:kernel_far_eps_proof} is proved, the restriction
\(x_1,x_3\in\mathbf C_\varepsilon\) can be removed. Indeed, by
\eqref{eq:kernel_bound_det},
\[
        |K_N(x_1,x_3)|
        \le
        \frac{C}{\sqrt{\det\Gamma_{x_1,x_3}}},
\]
and the latter function is integrable over \(\mathbf C^2\). Hence
\[
\begin{aligned}
        &\sup_{N\ge1}
        \int_{\mathbf C^2}
        \mathbf 1_{\{\|x_1-x_3\|>\delta\}}
        \mathbf 1_{\{x_1\notin\mathbf C_\varepsilon
        \ \text{or}\ x_3\notin\mathbf C_\varepsilon\}}
        |K_N(x_1,x_3)|\,\mathrm{d}x_1\mathrm{d}x_3        \\
        &\qquad\le
        C
        \int_{\mathbf C^2}
        \mathbf 1_{\{x_1\notin\mathbf C_\varepsilon
        \ \text{or}\ x_3\notin\mathbf C_\varepsilon\}}
        \frac{\mathrm{d}x_1\mathrm{d}x_3}{\sqrt{\det\Gamma_{x_1,x_3}}}
        \longrightarrow0
\end{aligned}
\]
as \(\varepsilon\downarrow0\), by the Dominated Convergence Theorem and Equation  \ref{eq:bound_frac_1_sqrt(det)}.

It remains to prove \eqref{eq:kernel_far_eps_proof}. Put
\[
        X_1=X(x_1),
        \qquad
        X_3=X(x_3),
        \qquad
        Z=(X_1,X_3)^\top .
\]
Conditionally on \(Z\), the vector
\[
        \mathbf U
        =
        \left(
        U^{(1)}_{x_1,x_2},
        U^{(2)}_{x_1,x_2},
        U^{(1)}_{x_3,x_4},
        U^{(2)}_{x_3,x_4}
        \right)^\top
\]
is Gaussian. Let
\[
        \Gamma_{x_1,x_3}
        =
        \operatorname{Cov}(Z).
\]
On the set
\[
        x_1,x_3\in\mathbf C_\varepsilon,
        \qquad
        \|x_1-x_3\|>\delta,
\]
the matrix \(\Gamma_{x_1,x_3}\) is uniformly non-degenerate: there exists
\(c_{\varepsilon,\delta}>0\) such that
\[
        \det\Gamma_{x_1,x_3}
        \ge c_{\varepsilon,\delta}.
\]

Let
\[
        \Lambda_{x_1,x_2,x_3,x_4}
        =
        \operatorname{Cov}(\mathbf U,Z).
\]
By Gaussian regression,
\[
        \mathbf U\mid Z
        \sim
        \mathcal N_4
        \left(
        \Lambda_{x_1,x_2,x_3,x_4}
        \Gamma_{x_1,x_3}^{-1}Z,
        \Sigma^{U}_{x_1,x_2,x_3,x_4}
        -
        \Lambda_{x_1,x_2,x_3,x_4}
        \Gamma_{x_1,x_3}^{-1}
        \Lambda_{x_1,x_2,x_3,x_4}^{\top}
        \right),
\]
where \(\Sigma^{U}_{x_1,x_2,x_3,x_4}\) denotes the covariance matrix of
\(\mathbf U\). By Lemma~\ref{Le:bound:corr}, which controls the correlation between two normalized increments, and by the following elementary bounds on the correlations between a normalized increment and a point evaluation of the
field,
\[
        |\rho_{x_1,x_2}|
        \le
        C_{\varepsilon}
        \left(
        \|x_2-x_1\|^{H}
        +
        \|x_2-x_1\|^{1-H}
        \right),
\]
\[
        |\rho_{x_3,x_4}|
        \le
        C_{\varepsilon}
        \left(
        \|x_4-x_3\|^{H}
        +
        \|x_4-x_3\|^{1-H}
        \right),
\]
and
\[
        |\nu_{x_1,x_2,x_3}|
        \le
        C_{\varepsilon,\delta}
        \|x_2-x_1\|^{1-H},
        \qquad
        |\nu_{x_3,x_4,x_1}|
        \le
        C_{\varepsilon,\delta}
        \|x_4-x_3\|^{1-H}
\]
uniformly for
\[
        x_1,x_3\in\mathbf C_\varepsilon,
        \qquad
        \|x_1-x_3\|>\delta,
\]
we have
\[
        \Sigma^{U}_{x_1,x_2,x_3,x_4}\to I_4,
        \qquad
        \Lambda_{x_1,x_2,x_3,x_4}\to0,
\]
as
\[
        \|x_2-x_1\|\vee\|x_4-x_3\|\to0.
\]
Indeed, these estimates follow directly from the covariance formula
\[
        \operatorname{Cov}(W(x),W(y))
        =
        \frac{1}{2}
        \left(
        \|x\|^{2H}+\|y\|^{2H}-\|x-y\|^{2H}
        \right),
\]
and from the fact that \(x\mapsto \|x\|^{2H}\) is \(C^1\) on every compact
set away from the origin. For instance, if \(x_1\in\mathbf C_\varepsilon\),
then, with \(h=x_2-x_1\),
\[
\begin{aligned}
        |\rho_{x_1,x_2}|
        &=
        \frac{1}{2}
        \left|
        \frac{
        -\|x_1+h\|^{2H}+\|x_1\|^{2H}+\|h\|^{2H}
        }
        {\|h\|^{H}\|x_1\|^{H}}
        \right|                                      \\
        &\le
        C_\varepsilon
        \frac{
        \|h\|+\|h\|^{2H}
        }
        {\|h\|^{H}}
        =
        C_\varepsilon
        \left(
        \|h\|^{1-H}
        +
        \|h\|^{H}
        \right).
\end{aligned}
\]
Similarly, if \(x_1,x_3\in\mathbf C_\varepsilon\) and
\(\|x_1-x_3\|>\delta\), then the map
\[
        h\mapsto
        \|x_1+h\|^{2H}-\|x_3-x_1-h\|^{2H}
\]
is \(C^1\) uniformly on the relevant compact set, which gives
\[
        |\nu_{x_1,x_2,x_3}|
        \le
        C_{\varepsilon,\delta}\|x_2-x_1\|^{1-H}.
\]

Consequently, if
\[
        X_1=N^{-H/2}z_1,
        \qquad
        X_3=N^{-H/2}z_3,
\]
with \(z_1,z_3\) fixed, then the conditional mean converges to \(0\) and the
conditional covariance converges to \(I_4\), uniformly on
\(\mathbf C_\varepsilon^2\cap\{\|x_1-x_3\|>\delta\}\), whenever the rescaled
edge lengths
\[
        \sqrt N\|x_2-x_1\|,
        \qquad
        \sqrt N\|x_4-x_3\|
\]
remain bounded. The contribution of unbounded rescaled edge lengths is
negligible by the exponential moment bound of Lemma~\ref{Le_Bound_Exp_R}.

We now use the local Delaunay scaling. Let
\[
        h=\sqrt N(x_2-x_1),
        \qquad
        g=\sqrt N(x_4-x_3).
\]
For \(x_1,x_3\in\mathbf C_\varepsilon\) with
\(\|x_1-x_3\|>\delta\), define the finite measure
\[
        Q^N_{x_1,x_3}(\mathrm{d}h,\mathrm{d}g)
        =
        p_{2,N}
        \left(
        x_1,x_1+\frac{h}{\sqrt N},
        x_3,x_3+\frac{g}{\sqrt N}
        \right)\,\mathrm{d}h\,\mathrm{d}g.
\]
The properties of the Poisson--Delaunay graph gives the vague convergence
\[
        Q^N_{x_1,x_3}
        \Longrightarrow
        Q\otimes Q,
\]
uniformly for
\[
        x_1,x_3\in\mathbf C_\varepsilon,
        \qquad
        \|x_1-x_3\|>\delta.
\]
Here \(Q\) is the oriented edge-intensity measure of the unit-intensity
Poisson--Delaunay graph. Its total mass is
\[
        Q(\mathbb R^2)=3,
\]
and the image of \(Q/3\) under \(h\mapsto\|h\|\) has density \(f_D\).

For \(z\in\mathbb R\) and \(h\neq0\), define
\[
        \Theta(z,h)
        =
        \int_{\mathbb R^2}
        \Psi
        \left(
        u,v,
        \frac{z}{\|h\|^{H}}
        \right)
        \frac{1}{2\pi}e^{-(u^2+v^2)/2}\,du\,dv .
\]
By the definition of \(F\),
\begin{equation}
        \int_{\mathbb R^2}\Theta(z,h)\,Q(\mathrm{d}h)
        =
        3F(z).
\label{eq:Theta_edge_average_kernel}
\end{equation}

For fixed \(z_1,z_3\in\mathbb R\), define
\[
\begin{aligned}
        \mathcal K_N(z_1,z_3;x_1,x_3)
        &=
        N^2
        \int_{(\mathbb R^2)^2}
        \mathbb E
        \left[
        r_{N;x_1,x_2}r_{N;x_3,x_4}
        \,\middle|\,
        X_1=N^{-H/2}z_1,\,
        X_3=N^{-H/2}z_3
        \right]                                      \\
        &\hspace{3cm}\times
        p_{2,N}(x_1,x_2,x_3,x_4)\,\mathrm{d}x_2\mathrm{d}x_4 .
\end{aligned}
\]
After the change of variables \(h=\sqrt N(x_2-x_1)\) and
\(g=\sqrt N(x_4-x_3)\), this becomes
\[
        \mathcal K_N(z_1,z_3;x_1,x_3)
        =
        \int_{\mathbb R^2\times\mathbb R^2}
        \mathcal R_N(z_1,z_3,h,g;x_1,x_3)
        Q^N_{x_1,x_3}(\mathrm{d}h,\mathrm{d}g),
\]
where \(\mathcal R_N\) denotes the conditional expectation of the product of
the two centered residuals under the conditioning
\[
        X_1=N^{-H/2}z_1,
        \qquad
        X_3=N^{-H/2}z_3.
\]
By the conditional Gaussian convergence above and
Lemma~\ref{Prop_Gen_Pod_Rosen}, for fixed \(z_1,z_3,h,g\),
\[
        \mathcal R_N(z_1,z_3,h,g;x_1,x_3)
        \longrightarrow
        \left\{\Theta(z_1,h)-F(z_1)\right\}
        \left\{\Theta(z_3,g)-F(z_3)\right\},
\]
uniformly for
\[
        x_1,x_3\in\mathbf C_\varepsilon,
        \qquad
        \|x_1-x_3\|>\delta.
\]
The bound \eqref{eq:Psi_basic_bound_kernel}, the decay of \(F\), and the
exponential tail of the typical Delaunay edge length give an integrable
dominating function in \((z_1,z_3,h,g)\). Hence dominated convergence and the
local Delaunay scaling imply
\[
\begin{aligned}
        \mathcal K_N(z_1,z_3;x_1,x_3)
        &\longrightarrow
        \int
        \left\{\Theta(z_1,h)-F(z_1)\right\}
        \left\{\Theta(z_3,g)-F(z_3)\right\}
        Q(\mathrm{d}h)Q(\mathrm{d}g)                                    \\
        &=
        \left[
        \int \Theta(z_1,h)\,Q(\mathrm{d}h)-3F(z_1)
        \right]
        \left[
        \int \Theta(z_3,g)\,Q(\mathrm{d}g)-3F(z_3)
        \right]                                      \\
        &=0,
\end{aligned}
\]
where the last equality follows from
\eqref{eq:Theta_edge_average_kernel}.

Finally, \(K_N(x_1,x_3)\) is obtained from \(\mathcal K_N\) by integrating
with respect to the density of
\[
        \left(
        N^{H/2}X_1,
        N^{H/2}X_3
        \right).
\]
On \(\mathbf C_\varepsilon^2\cap\{\|x_1-x_3\|>\delta\}\), these densities are
uniformly bounded by a Gaussian density with uniformly non-degenerate
covariance. The tails in \(z_1,z_3\) are controlled by the decay of \(F\) and
by the preceding domination. Therefore, using Lemma \ref{Prop_Gen_Pod_Rosen} gives
\[
        \sup_{\substack{x_1,x_3\in\mathbf C_\varepsilon\\
        \|x_1-x_3\|>\delta}}
        |K_N(x_1,x_3)|
        \longrightarrow0.
\]
This proves \eqref{eq:kernel_far_eps_proof}. Letting
\(\varepsilon\downarrow0\) using the determinant domination established at
the beginning of the proof yields
\[
        \int_{\mathbf C^2}
        \mathbf 1_{\{\|x_1-x_3\|>\delta\}}
        K_N(x_1,x_3)\,\mathrm{d}x_1\mathrm{d}x_3
        \longrightarrow0.
\]
This is \eqref{eq:kernel_bound_far_integral}, and the proof of the lemma is
complete.
\end{prooft}

\subsection{Proof of Proposition \protect\ref{prop:Brownian_parts}}
\label{subsec:proof_negligible_gaussian_parts}

We prove that
\[
        N^{-(1-H)/2}
        \left(
        V^{(1)}_{2,N}+V^{(2)}_{2,N}
        \right)
        \xrightarrow{L^2}0,
        \qquad H\in(0,1/2),
\]
or equivalently
\[
        N^{-(1-H)/2}V^{(i)}_{2,N}
        \xrightarrow{L^2}0,
        \qquad i=1,2.
\]

We shall use the following auxiliary estimate, proved in
Section~\ref{sec:variance_truncated_quadratic_variation} of the Supplementary Material.

\begin{proposition}[Variance bound for the truncated quadratic variation]
\label{prop:variance_truncated_quadratic_variation}
Let \(P_1\) be a homogeneous Poisson point process with intensity \(1\) on
\(\mathbb R^2\), and set
\[
        \mathbf C_N=\left(-\frac{\sqrt N}{2},\frac{\sqrt N}{2}\right]^2 .
\]
Let \(E'_N\) denote the set of oriented Delaunay edges of the graph generated
by \(P_1\cap \mathbf C_N\). Let \(W\) and \(V\) be two independent and identically
distributed isotropic fractional Brownian fields with Hurst parameter
\(H\in(0,1/2)\) and scale parameter \(\sigma^2\). For an oriented edge
\(e=(x,y)\in E'_N\), define
\[
        U^W_e
        =
        \frac{W(y)-W(x)}{\sigma \|y-x\|^H}.
\]
For \(\varepsilon\in\{-1,+1\}\), define the sign mark
\[
        M_\varepsilon(x)
        =
        \mathbf 1_{\{\varepsilon (V(x)-W(x))>0\}},
\]
and the truncated quadratic-variation statistic
\[
        T^{\varepsilon}_{2,N}
        =
        \frac{1}{\sqrt{|E'_N|}}
        \sum_{e=(x,y)\in E'_N}
        \left\{
        (U^W_e)^2-1
        \right\}
        M_\varepsilon(x).
\]
Then there exists a constant \(C<\infty\), independent of \(N\), such that
\[
        \mathbb E\left[
        \left(T^{\varepsilon}_{2,N}\right)^2
        \right]
        \le C N^{1-2H},
        \qquad N\ge 1,\quad \varepsilon\in\{-1,+1\}.
\]
Consequently,
\[
        N^{-(1-H)/2}T^{\varepsilon}_{2,N}
        \xrightarrow{L^2}0 .
\]
\end{proposition}

We now deduce Proposition~\ref{prop:Brownian_parts} from
Proposition~\ref{prop:variance_truncated_quadratic_variation}.

Let \(P_1\) be a homogeneous Poisson point process with intensity \(1\) on
\(\mathbb R^2\). By the scaling property of Poisson point processes,
$P_N \stackrel{d}{=} N^{-1/2}P_1$. Moreover, the Delaunay graph is invariant under
similarities. Hence the oriented Delaunay edges in the fixed window
\(\mathbf C\) may be represented as
\[
        E_N
        \stackrel{d}{=}
        \left\{
        \left(\frac{x}{\sqrt N},\frac{y}{\sqrt N}\right)
        :
        (x,y)\in E'_N
        \right\}.
\]

Let us first consider \(V^{(1)}_{2,N}\). Recall that
\[
        V^{(1)}_{2,N}
        =
        \frac{1}{\sqrt{|E_N|}}
        \sum_{(x_1,x_2)\in E_N}
        \left\{
        \left(U^{(1)}_{x_1,x_2}\right)^2-1
        \right\}
        \mathbf 1_{\{W^{(2)}(x_1)-W^{(1)}(x_1)<0\}},
\]
where
\[
        U^{(1)}_{x_1,x_2}
        =
        \frac{W^{(1)}(x_2)-W^{(1)}(x_1)}
        {\sigma \|x_2-x_1\|^H}.
\]
For \(x\in C_N\), define the rescaled fields
\[
        \widetilde W^{(i)}_N(x)
        =
        N^{H/2}W^{(i)}\left(\frac{x}{\sqrt N}\right),
        \qquad i=1,2.
\]
By the self-similarity of isotropic fractional Brownian fields and by the
stationarity of their increments,
\[
        \left(\widetilde W^{(1)}_N,\widetilde W^{(2)}_N\right)
        \stackrel{d}{=}
        \left(W^{(1)},W^{(2)}\right).
\]
Furthermore, for every edge \((x,y)\in E'_N\),
\[
        \frac{
        W^{(1)}(y/\sqrt N)-W^{(1)}(x/\sqrt N)
        }
        {\sigma \|y/\sqrt N-x/\sqrt N\|^H}
        =
        \frac{
        \widetilde W^{(1)}_N(y)-\widetilde W^{(1)}_N(x)
        }
        {\sigma \|y-x\|^H}.
\]
The sign mark is also preserved by the same rescaling, since
\[
        W^{(2)}(x/\sqrt N)-W^{(1)}(x/\sqrt N)<0
\]
is equivalent to
\[
        \widetilde W^{(2)}_N(x)-\widetilde W^{(1)}_N(x)<0.
\]
Therefore,
\[
        V^{(1)}_{2,N}
        \stackrel{d}{=}
        T^{-}_{2,N},
\]
with \(W=\widetilde W^{(1)}_N\), \(V=\widetilde W^{(2)}_N\), and
\(\varepsilon=-1\) in Proposition
\ref{prop:variance_truncated_quadratic_variation}. Hence
\[
        \mathbb E\left[
        \left(V^{(1)}_{2,N}\right)^2
        \right]
        \le C N^{1-2H}.
\]
It follows that
\[
        \mathbb E\left[
        \left(
        N^{-(1-H)/2}V^{(1)}_{2,N}
        \right)^2
        \right]
        \le
        C N^{-(1-H)}N^{1-2H}
        =
        C N^{-H}
        \longrightarrow 0.
\]
Thus
\[
        N^{-(1-H)/2}V^{(1)}_{2,N}
        \xrightarrow{L^2}0.
\]

The proof for \(V^{(2)}_{2,N}\) is identical. Indeed,
\[
        V^{(2)}_{2,N}
        =
        \frac{1}{\sqrt{|E_N|}}
        \sum_{(x_1,x_2)\in E_N}
        \left\{
        \left(U^{(2)}_{x_1,x_2}\right)^2-1
        \right\}
        \mathbf 1_{\{W^{(2)}(x_1)-W^{(1)}(x_1)>0\}},
\]
and the same Poisson rescaling and fractional Brownian self-similarity give
\[
        V^{(2)}_{2,N}
        \stackrel{d}{=}
        T^{+}_{2,N},
\]
possibly after exchanging the two independent copies. Therefore
\[
        \mathbb E\left[
        \left(V^{(2)}_{2,N}\right)^2
        \right]
        \le C N^{1-2H},
\]
and consequently
\[
        N^{-(1-H)/2}V^{(2)}_{2,N}
        \xrightarrow{L^2}0.
\]

Finally, by the elementary inequality \((a+b)^2\le 2a^2+2b^2\),
\[
\begin{aligned}
        \mathbb E\left[
        \left(
        N^{-(1-H)/2}
        \left(
        V^{(1)}_{2,N}+V^{(2)}_{2,N}
        \right)
        \right)^2
        \right]
        &\le
        2N^{-(1-H)}
        \mathbb E\left[\left(V^{(1)}_{2,N}\right)^2\right]  \\
        &\quad+
        2N^{-(1-H)}
        \mathbb E\left[\left(V^{(2)}_{2,N}\right)^2\right]  \\
        &\le
        C N^{-H}.
\end{aligned}
\]
Since \(H>0\), the right-hand side converges to \(0\). Hence
\[
        N^{-(1-H)/2}
        \left(
        V^{(1)}_{2,N}+V^{(2)}_{2,N}
        \right)
        \xrightarrow{L^2}0.
\]
This proves
Proposition~\ref{prop:Brownian_parts}.

\section{Conclusion}
\label{sec:conclusion}

This paper establishes a limit theorem for squared increment sums of the
pointwise maximum of two independent isotropic fractional Brownian fields
observed in a fixed two-dimensional domain. The observation design is random
and is modeled by a homogeneous Poisson point process, while increments are
computed along oriented edges of the associated Delaunay triangulation. The
main conclusion is that the asymptotic behavior of the quadratic variation of
the maximum field differs fundamentally from the Gaussian fluctuation regime
obtained for a single fractional Brownian field.

The key point is the decomposition of the normalized squared increment sum into
two Gaussian-type terms and one transition term. The Gaussian-type terms
correspond to edges for which the same fractional Brownian field realizes the
maximum at the initial endpoint. At the scale relevant for the maximum field,
these contributions are negligible. By contrast, the transition term collects
the contribution of edges for which the identity of the maximizer may change
along the edge. This term is non-zero only when the difference field
\(
        W^{(2\backslash 1)} = W^{(2)} - W^{(1)}
\)
is of the order of the local increment scale. Since a typical
Poisson--Delaunay edge in the fixed window has length of order \(N^{-1/2}\),
this mechanism localizes the leading contribution in a shrinking band of width
\(N^{-H/2}\) around the random interface
\(
        \{x\in \mathbf{C} : W^{(1)}(x)=W^{(2)}(x)\}.
\)

The main theorem shows that, for \(H\in(0,1/2)\),
\[
        \frac{\sqrt{3}}{3}\,
        N^{-(1-H)/2}
        V^{(W^\vee)}_{2,N}
        \xrightarrow[L^2]{}
        c_{V_2} L_{W^{(2\backslash 1)}}(0).
\]
Thus the limiting object is not a centered Gaussian random variable, but the
local time at zero of the difference field. This local time measures, in the
occupation-density sense, the size of the random contact set where the two
fractional Brownian fields coincide. The result therefore identifies a
genuinely geometric asymptotic regime: the leading contribution is determined
by the random interface between the two competing fields, rather than by the
accumulation of weakly dependent Gaussian increment fluctuations.

Several technical ingredients are combined to prove this result. On the
stochastic-geometry side, the Poisson--Delaunay structure provides the
edge-intensity normalization, the distribution of the typical edge length, and
tractable expectation formulas through the Slivnyak--Mecke formula. These tools
make it possible to replace local sums over Delaunay edges by averaged
quantities involving the typical edge. The proof also uses stabilization
properties of Poisson--Delaunay functionals and bounds on the probability that
prescribed pairs of points form Delaunay edges.

On the Gaussian side, the main tools are occupation-density arguments and an
\(L^2\)-Fourier representation of the local time of the difference field. The
proof of the transition limit relies on a residual covariance analysis: after
subtracting the averaged edge contribution, one controls a covariance kernel by
a determinant bound for the bivariate Gaussian vector
\(
        \bigl(
        W^{(2\backslash 1)}(x_1),
        W^{(2\backslash 1)}(x_3)
        \bigr).
\)
The near-diagonal part is handled through local nondeterminism-type estimates,
while the far-field part is treated by conditional Gaussian regression, local
Delaunay scaling, and a Gaussian comparison argument.

Finally, the negligibility of the Gaussian-type terms is obtained through a
variance bound for a truncated quadratic variation in the rescaled Poisson
framework. This step uses the self-similarity and stationary increments of the
fractional Brownian field, the summability of squared increment correlations
for \(H<1/2\), and regression estimates controlling the dependence between
normalized increments and the sign of the difference field. Together, these
arguments show that the pointwise maximum of two fractional Brownian fields
induces a new fixed-domain asymptotic regime for Delaunay-edge quadratic
variations, governed by local time and by the geometry of the contact set.

\bibliographystyle{abbrv}
\bibliography{biblio_infill}

@article {Chan&Wood00,
    AUTHOR = {Chan, Grace and Wood, Andrew T. A.},
     TITLE = {Increment-based estimators of fractal dimension for
              two-dimensional surface data},
   JOURNAL = {Statist. Sinica},
  FJOURNAL = {Statistica Sinica},
    VOLUME = {10},
      YEAR = {2000},
    NUMBER = {2},
     PAGES = {343--376},
      ISSN = {1017-0405,1996-8507},
   MRCLASS = {62M40 (60G18 60G60 62F12)},
  MRNUMBER = {1769748},
}

@article {Chenavier&Robert25a,
    AUTHOR = {Chenavier, Nicolas and Robert, Christian Y.},
     TITLE = {Central limit theorems for squared increment sums of fractional {B}rownian fields based on a {D}elaunay triangulation in {$2D$}},
   JOURNAL = {WP},
  FJOURNAL = {Working paper},
      YEAR = {2026},
}

@article {Chenavier&Robert25c,
    AUTHOR = {Chenavier, Nicolas and Robert, Christian Y.},
     TITLE = {Asymptotic properties of maximum composite likelihood estimators for max-stable {B}rown-{R}esnick random fields over a fixed-domain},
   JOURNAL = {WP},
  FJOURNAL = {Working paper},
      YEAR = {2026},
}

@article {Geman&Horowitz80,
    AUTHOR = {Geman, Donald and Horowitz, Joseph},
     TITLE = {Occupation densities},
   JOURNAL = {Ann. Probab.},
  FJOURNAL = {The Annals of Probability},
    VOLUME = {8},
      YEAR = {1980},
    NUMBER = {1},
     PAGES = {1--67},
      ISSN = {0091-1798,2168-894X},
   MRCLASS = {60J55 (26A27 60G15 60G17)},
  MRNUMBER = {556414},
MRREVIEWER = {Simeon\ M.\ Berman},
       URL =
              {http://links.jstor.org/sici?sici=0091-1798(198002)8:1<1:OD>2.0.CO;2-M&origin=MSN},
}

@article {Jaramillo21,
    AUTHOR = {Jaramillo, Arturo and Nourdin, Ivan and Peccati, Giovanni},
     TITLE = {Approximation of fractional local times: zero energy and
              derivatives},
   JOURNAL = {Ann. Appl. Probab.},
  FJOURNAL = {The Annals of Applied Probability},
    VOLUME = {31},
      YEAR = {2021},
    NUMBER = {5},
     PAGES = {2143--2191},
      ISSN = {1050-5164,2168-8737},
   MRCLASS = {60G22 (60F17 60H07 60J55)},
  MRNUMBER = {4332693},
       DOI = {10.1214/20-aap1643},
       URL = {https://doi.org/10.1214/20-aap1643},
}

@article {Podolskij&Rosenbaum18,
    AUTHOR = {Podolskij, Mark and Rosenbaum, Mathieu},
     TITLE = {Comment on: Limit of random
measures associated with the increments of a {B}rownian semimartingale.
Asymptotic behavior of local times related statistics for fractional
{B}rownian motion},
   JOURNAL = {Journal of Financial Econometrics},
    VOLUME = {16},
      YEAR = {2018},
    NUMBER = {4},
     PAGES = {588--598},
}

@article {Robert20,
    AUTHOR = {Robert, Christian Y.},
     TITLE = {Power variations for a class of {B}rown-{R}esnick processes},
   JOURNAL = {Extremes},
  FJOURNAL = {Extremes. Statistical Theory and Applications in Science,
              Engineering and Economics},
    VOLUME = {23},
      YEAR = {2020},
    NUMBER = {2},
     PAGES = {215--244},
      ISSN = {1386-1999,1572-915X},
   MRCLASS = {60G70 (60F05 60G17)},
  MRNUMBER = {4102383},
MRREVIEWER = {Peter\ Kern},
       DOI = {10.1007/s10687-020-00373-4},
       URL = {https://doi.org/10.1007/s10687-020-00373-4},
}

@book {Schneider&Weil08,
    AUTHOR = {Schneider, Rolf and Weil, Wolfgang},
     TITLE = {Stochastic and integral geometry},
    SERIES = {Probability and its Applications (New York)},
 PUBLISHER = {Springer-Verlag, Berlin},
      YEAR = {2008},
     PAGES = {xii+693},
      ISBN = {978-3-540-78858-4},
   MRCLASS = {60-02 (52A22 60D05 60G55 62M30)},
  MRNUMBER = {2455326},
MRREVIEWER = {V.\ K.\ Ohanyan},
       DOI = {10.1007/978-3-540-78859-1},
       URL = {https://doi.org/10.1007/978-3-540-78859-1},
}

@article {Zhu&Stein02,
    AUTHOR = {Zhu, Zhengyuan and Stein, Michael L.},
     TITLE = {Parameter estimation for fractional {B}rownian surfaces},
   JOURNAL = {Statist. Sinica},
  FJOURNAL = {Statistica Sinica},
    VOLUME = {12},
      YEAR = {2002},
    NUMBER = {3},
     PAGES = {863--883},
      ISSN = {1017-0405,1996-8507},
   MRCLASS = {28A80 (60G60)},
  MRNUMBER = {1929968},
}

\newpage

\appendix

\noindent
{\LARGE\bf Supplementary material}

\section{Existence of the \(L^2\) Fourier representation of the local time}
\label{app:L2_local_time}
\label{sec:L2_local_time}

In this section we provide a self-contained proof of the \(L^2\)
Fourier representation of the local time at level \(0\) for the Gaussian
field \(W^{(2\backslash 1)}\). The same argument applies to any fixed level
\(\ell\in\mathbb R\), up to the additional phase factor \(e^{-i\xi\ell}\) in
the Fourier representation.

Let
\[
        X(x)=W^{(2\backslash 1)}(x)
        =
        W^{(2)}(x)-W^{(1)}(x),
        \qquad x\in\mathbb R^2.
\]
Then \(X\) is a centered Gaussian field with stationary increments and
covariance
\[
        \operatorname{Cov}(X(x),X(y))
        =
        \sigma^2
        \left(
        \|x\|^{2H}+\|y\|^{2H}-\|x-y\|^{2H}
        \right).
\]
In particular,
\begin{equation}
        \Delta(x,y)
        :=
        \operatorname{Var}(X(x)-X(y))
        =
        2\sigma^2\|x-y\|^{2H},
        \qquad x,y\in\mathbb R^2.
\label{eq:Delta_app}
\end{equation}

\begin{proposition}
\label{prop:L2_Fourier_local_time}
There exists a random variable \(L_X(0)\in L^2\) such that, for
\(\varepsilon>0\),
\[
        L_\varepsilon
        :=
        \int_{\mathbf C}
        \varphi_\varepsilon(X(x))\,\mathrm{d}x,
        \qquad
        \varphi_\varepsilon(u)
        :=
        (2\pi\varepsilon)^{-1/2}
        \exp\left(-\frac{u^2}{2\varepsilon}\right),
\]
converges to \(L_X(0)\) in \(L^2\) as
\(\varepsilon\downarrow0\). Moreover,
\begin{equation}
        L_X(0)
        =
        \frac{1}{2\pi}
        \lim_{M\to\infty}
        \int_{-M}^{M}
        \int_{\mathbf C}
        \exp\left(i\xi X(x)\right)\,\mathrm{d}x\,\mathrm{d}\xi,
        \qquad \text{in }L^2.
\label{eq:Fourier_L2_app}
\end{equation}
\end{proposition}

We first record a determinant estimate for the two-dimensional Gaussian vector
\((X(x),X(y))\). It is a two-point form of local nondeterminism and is the
key integrability estimate used below.

\begin{lemma}
\label{lem:det_lower_bound_app}
There exists a constant \(c>0\), depending only on \(H\), \(\sigma\) and
\(\mathbf C\), such that, for all \(x,y\in\mathbf C\),
\begin{equation}
        \det\operatorname{Cov}(X(x),X(y))
        \ge
        c\,
        \bigl(\|x\|\wedge\|y\|\bigr)^{2H}
        \|x-y\|^{2H} .
\label{eq:det_lower_bound_app}
\end{equation}
Consequently,
\begin{equation}
        \int_{\mathbf C^2}
        \frac{\mathrm{d}x\,\mathrm{d}y}
        {\sqrt{\det\operatorname{Cov}(X(x),X(y))}}
        <\infty .
\label{eq:det_integrability_app}
\end{equation}
\end{lemma}

\begin{proof}
We first prove the lower bound. For a centered Gaussian vector
\((Z_1,Z_2)\), one has
\[
        \det\operatorname{Cov}(Z_1,Z_2)
        =
        \operatorname{Var}(Z_2)
        \operatorname{Var}(Z_1\mid Z_2).
\]
Applied to \(Z_1=X(x)\), \(Z_2=X(y)\), this gives
\[
        \det\operatorname{Cov}(X(x),X(y))
        =
        \operatorname{Var}(X(y))
        \operatorname{Var}(X(x)\mid X(y)).
\]
The isotropic fractional Brownian field satisfies the two-point local
nondeterminism bound
\[
        \operatorname{Var}(X(x)\mid X(y))
        \ge
        c_1
        \bigl(
        \|x\|^{2H}\wedge\|x-y\|^{2H}
        \bigr),
\]
for all \(x,y\in\mathbf C\). This estimate follows directly from the explicit
covariance formula by homogeneity and a compactness argument: after
normalizing by \(\|x\|\), the ratio
\[
        \frac{\operatorname{Var}(X(x)\mid X(y))}
        {\|x\|^{2H}\wedge\|x-y\|^{2H}}
\]
extends continuously and remains strictly positive away from the degenerate
cases \(x=0\) and \(x=y\), where the same lower bound follows from the
asymptotics of the variogram \(\|x-y\|^{2H}\). Hence
\[
        \det\operatorname{Cov}(X(x),X(y))
        \ge
        c_2
        \|y\|^{2H}
        \bigl(
        \|x\|^{2H}\wedge\|x-y\|^{2H}
        \bigr).
\]
By symmetry, the same argument with \(x\) and \(y\) exchanged gives
\[
        \det\operatorname{Cov}(X(x),X(y))
        \ge
        c_2
        \|x\|^{2H}
        \bigl(
        \|y\|^{2H}\wedge\|x-y\|^{2H}
        \bigr).
\]
Combining the two estimates yields
\[
        \det\operatorname{Cov}(X(x),X(y))
        \ge
        c\,
        \bigl(\|x\|\wedge\|y\|\bigr)^{2H}
        \|x-y\|^{2H},
\]
possibly after reducing \(c\), since \(x,y\) belong to the bounded set
\(\mathbf C\).

It remains to prove integrability. Let \(p=H\in(0,1)\). From
\eqref{eq:det_lower_bound_app},
\[
        \frac{1}
        {\sqrt{\det\operatorname{Cov}(X(x),X(y))}}
        \le
        C
        \bigl(\|x\|\wedge\|y\|\bigr)^{-p}
        \|x-y\|^{-p}.
\]
Using
\[
        \bigl(\|x\|\wedge\|y\|\bigr)^{-p}
        \le
        \|x\|^{-p}+\|y\|^{-p},
\]
we get
\[
\begin{aligned}
        \int_{\mathbf C^2}
        \frac{\mathrm{d}x\,\mathrm{d}y}
        {\sqrt{\det\operatorname{Cov}(X(x),X(y))}}
        &\le
        C
        \int_{\mathbf C^2}
        \left(\|x\|^{-p}+\|y\|^{-p}\right)
        \|x-y\|^{-p}\,\mathrm{d}x\,\mathrm{d}y        \\
        &=
        2C
        \int_{\mathbf C}
        \|x\|^{-p}
        \left(
        \int_{\mathbf C}
        \|x-y\|^{-p}\,\mathrm{d}y
        \right)\mathrm{d}x.
\end{aligned}
\]
Since \(p<2\), the inner integral is uniformly bounded in \(x\in\mathbf C\),
and since again \(p<2\),
\[
        \int_{\mathbf C}\|x\|^{-p}\,\mathrm{d}x<\infty.
\]
This proves \eqref{eq:det_integrability_app}.
\end{proof}

\strut

\begin{proof}[Proof of Proposition~\ref{prop:L2_Fourier_local_time}]
For \(\varepsilon>0\), the Gaussian kernel admits the Fourier representation
\[
        \varphi_\varepsilon(u)
        =
        \frac{1}{2\pi}
        \int_{\mathbb R}
        e^{i\xi u}e^{-\varepsilon\xi^2/2}\,\mathrm{d}\xi .
\]
Therefore, by Fubini,
\begin{equation}
        L_\varepsilon
        =
        \frac{1}{2\pi}
        \int_{\mathbb R}
        \int_{\mathbf C}
        e^{i\xi X(x)}
        e^{-\varepsilon\xi^2/2}
        \,\mathrm{d}x\,\mathrm{d}\xi .
\label{eq:L_eps_Fourier_app}
\end{equation}

We first show that \((L_\varepsilon)_{\varepsilon>0}\) is Cauchy in
\(L^2\). Let \(\varepsilon,\eta>0\), and set
\[
        K_{\varepsilon,\eta}(u)
        =
        e^{-\varepsilon u^2/2}
        -
        e^{-\eta u^2/2}.
\]
Using \eqref{eq:L_eps_Fourier_app} and Tonelli's theorem, we obtain
\begin{equation}
        \mathbb E\left[
        (L_\varepsilon-L_\eta)^2
        \right]
        =
        \frac{1}{(2\pi)^2}
        \int_{\mathbb R^2}
        \int_{\mathbf C^2}
        K_{\varepsilon,\eta}(\xi)
        K_{\varepsilon,\eta}(\xi')     
        \mathbb E\left[
        \exp\left(
        i\xi X(x)+i\xi'X(y)
        \right)
        \right]
        \,\mathrm{d}x\,\mathrm{d}y\,\mathrm{d}\xi\,\mathrm{d}\xi' .
\label{eq:Cauchy_expansion_app}
\end{equation}
For fixed \(x,y\in\mathbf C\), let \(\Sigma_{x,y}\) be the covariance matrix
of \((X(x),X(y))\). Since this vector is centered Gaussian,
\[
        \mathbb E\left[
        \exp\left(
        i\xi X(x)+i\xi'X(y)
        \right)
        \right]
        =
        \exp\left(
        -\frac12
        (\xi,\xi')\Sigma_{x,y}(\xi,\xi')^\top
        \right).
\]
Moreover,
\[
        |K_{\varepsilon,\eta}(u)|\le 1,
        \qquad
        K_{\varepsilon,\eta}(u)\to0
\]
pointwise as \((\varepsilon,\eta)\downarrow(0,0)\). Hence, by dominated
convergence, it is enough to prove that
\begin{equation}
        \int_{\mathbb R^2}
        \int_{\mathbf C^2}
        \exp\left(
        -\frac12
        (\xi,\xi')\Sigma_{x,y}(\xi,\xi')^\top
        \right)
        \,\mathrm{d}x\,\mathrm{d}y\,\mathrm{d}\xi\,\mathrm{d}\xi'
        <\infty .
\label{eq:dom_need_app}
\end{equation}
For any positive definite \(2\times2\) matrix \(\Sigma\),
\[
        \int_{\mathbb R^2}
        \exp\left(
        -\frac12 z^\top\Sigma z
        \right)\,\mathrm{d}z
        =
        \frac{2\pi}{\sqrt{\det\Sigma}}.
\]
Thus \eqref{eq:dom_need_app} follows directly from
\eqref{eq:det_integrability_app}. Consequently,
\((L_\varepsilon)_{\varepsilon>0}\) is Cauchy in \(L^2\). We denote
its limit by \(L_X(0)\).

We now identify this limit with the truncated Fourier representation. For
\(M>0\), define
\[
        L^{(M)}
        :=
        \frac{1}{2\pi}
        \int_{-M}^{M}
        \int_{\mathbf C}
        e^{i\xi X(x)}
        \,\mathrm{d}x\,\mathrm{d}\xi .
\]
Repeating the preceding argument with
\[
        \mathbf 1_{[-M,M]}-\mathbf 1_{[-M',M']}
\]
in place of \(K_{\varepsilon,\eta}\), and using the same dominating function,
shows that \((L^{(M)})_{M>0}\) is Cauchy in \(L^2\). Let
\(\widetilde L\) denote its \(L^2\)-limit.

Finally, for fixed \(M\),
\[
        \frac{1}{2\pi}
        \int_{-M}^{M}
        \int_{\mathbf C}
        e^{i\xi X(x)}
        e^{-\varepsilon\xi^2/2}
        \,\mathrm{d}x\,\mathrm{d}\xi
        \xrightarrow[\varepsilon\downarrow0]{L^2}
        L^{(M)}.
\]
On the other hand, letting \(M\to\infty\) in the damped representation
recovers \(L_\varepsilon\), again by the domination established above. The two
limits are therefore compatible, and the \(L^2\)-limit of \(L^{(M)}\) must be
the same as the \(L^2\)-limit of \(L_\varepsilon\). Hence
\[
        \widetilde L=L_X(0).
\]
This proves
\[
        L_X(0)
        =
        \frac{1}{2\pi}
        \lim_{M\to\infty}
        \int_{-M}^{M}
        \int_{\mathbf C}
        e^{i\xi X(x)}
        \,\mathrm{d}x\,\mathrm{d}\xi
\]
in \(L^2\), which is \eqref{eq:Fourier_L2_app}.
\end{proof}

\section{Proof of Proposition~\ref{prop:variance_truncated_quadratic_variation}}
\label{sec:variance_truncated_quadratic_variation}

We work in the
rescaled Poisson framework. Let \(P_1\) be a homogeneous Poisson point process
with intensity \(1\) on \(\mathbb R^2\), and set
\[
       \mathbf C_N=\left(-\frac{\sqrt N}{2},\frac{\sqrt N}{2}\right]^2 .
\]
Let \(E'_N\) denote the set of oriented Delaunay edges of the graph generated
by \(P_1\cap \mathbf C_N\) (more precisely with leftmost point in \(P_1\cap \mathbf C_N\)). For an oriented edge \(e=(x,y)\in E'_N\), write
\[
        \ell(e)=\|y-x\|.
\]
Let \(W\) and \(V\) be two independent isotropic fractional Brownian fields
with common Hurst parameter \(H\in(0,1/2)\) and common scale parameter
\(\sigma^2\). Set
\[
        D(x)=V(x)-W(x),
        \qquad x\in\mathbb R^2 .
\]
For \(\varepsilon\in\{-1,+1\}\), define
\[
        M_\varepsilon(x)
        =
        \mathbf 1_{\{\varepsilon D(x)>0\}},
\]
and, for every oriented edge \(e=(x,y)\in E'_N\),
\[
        U_e
        =
        \frac{W(y)-W(x)}{\sigma \ell(e)^H},
        \qquad
        Y^\varepsilon_e
        =
        \left(U_e^2-1\right)M_\varepsilon(x).
\]
The truncated quadratic-variation statistic is
\[
        T^\varepsilon_{2,N}
        =
        \frac{1}{\sqrt{|E'_N|}}
        \sum_{e\in E'_N}Y^\varepsilon_e .
\]
We prove that there exists a constant \(C<\infty\), independent of \(N\), such
that
\begin{equation}
        \mathbb E\left[
        \left(T^\varepsilon_{2,N}\right)^2
        \right]
        \le
        C N^{1-2H},
        \qquad N\ge1 .
\label{eq:variance_bound_truncated_goal}
\end{equation}

\subsection{Covariance decomposition}

All conditional expectations and covariances below are taken with respect to
the Gaussian fields, conditionally on the Poisson point process \(P_1\). For
fixed \(P_1\), the edge set \(E'_N\) is deterministic.

We first record a cancellation. For every fixed edge \(e=(x,y)\),
\[
        \mathbb E\left[
        Y^\varepsilon_e\mid P_1
        \right]=0.
\]
Indeed, under the transformation
\[
        (W,V)\mapsto (-W,-V),
\]
the joint distribution of \((W,V)\) is unchanged. Moreover \(U_e^2-1\) is
invariant, whereas \(M_\varepsilon(x)\) is transformed into
\(1-M_\varepsilon(x)\), up to the null event \(\{D(x)=0\}\). Hence
\[
        \mathbb E\left[
        M_\varepsilon(x)(U_e^2-1)\mid P_1
        \right]
        =
        \mathbb E\left[
        (1-M_\varepsilon(x))(U_e^2-1)\mid P_1
        \right].
\]
Adding the two equal quantities gives
\[
        2\mathbb E\left[
        M_\varepsilon(x)(U_e^2-1)\mid P_1
        \right]
        =
        \mathbb E\left[
        U_e^2-1\mid P_1
        \right]
        =
        0,
\]
because \(U_e\) is standard normal conditionally on \(P_1\). Consequently,
\begin{equation}
        \mathbb E\left[
        \left(T^\varepsilon_{2,N}\right)^2
        \,\middle|\,P_1
        \right]
        =
        \frac{1}{|E'_N|}
        \sum_{e,e'\in E'_N}
        \operatorname{Cov}
        \left(
        Y^\varepsilon_e,Y^\varepsilon_{e'}
        \,\middle|\,P_1
        \right).
\label{eq:second_moment_covariance_decomposition_start}
\end{equation}

Let
\[
        e=(x_1,x_2),
        \qquad
        e'=(x_3,x_4)
\]
be two oriented edges of \(E'_N\). We write
\[
        U=U_e,
        \qquad
        U'=U_{e'},
        \qquad
        M_1=M_\varepsilon(x_1),
        \qquad
        M_3=M_\varepsilon(x_3),
\]
and
\[
        \rho=\operatorname{Cov}(U,U'\mid P_1)
        =
        \operatorname{Corr}(U,U'\mid P_1).
\]
The equality between covariance and correlation follows from
\[
        \operatorname{Var}(U\mid P_1)
        =
        \operatorname{Var}(U'\mid P_1)
        =
        1.
\]

Introduce
\[
        Z=(W(x_1),W(x_3))^\top .
\]
Since \(V\) is independent of \(W\), and since \(U,U'\) are increments of
\(W\) only, the conditional distribution of \((U,U')\) given
\[
        \mathcal G
        =
        \sigma\{W(x_1),V(x_1),W(x_3),V(x_3)\}
\]
depends on \(\mathcal G\) only through \(Z\). Define the Gaussian regressions
\[
        m=\mathbb E[U\mid Z,P_1],
        \qquad
        m'=\mathbb E[U'\mid Z,P_1].
\]
Then
\[
        U=m+\eta,
        \qquad
        U'=m'+\eta',
\]
where \((\eta,\eta')\) is a centered Gaussian vector independent of \(Z\).
Let
\[
        c=\operatorname{Cov}(\eta,\eta'\mid P_1).
\]
Since \(U\) and \(U'\) have unit conditional variances,
\[
        c
        =
        \rho-\mathbb E[mm'\mid P_1].
\]

Set
\[
        \Delta
        =
        m^2-\mathbb E[m^2\mid P_1],
        \qquad
        \Delta'
        =
        (m')^2-\mathbb E[(m')^2\mid P_1],
\]
and
\[
        \Xi
        =
        mm'-\mathbb E[mm'\mid P_1].
\]
By the law of total covariance with respect to \(\mathcal G\),
\begin{equation*}
        \operatorname{Cov}
        \left(
        Y^\varepsilon_e,Y^\varepsilon_{e'}
        \,\middle|\,P_1
        \right)                                                   
        =
        \mathbb E\left[
        \operatorname{Cov}
        \left(
        Y^\varepsilon_e,Y^\varepsilon_{e'}
        \,\middle|\,\mathcal G,P_1
        \right)
        \,\middle|\,P_1
        \right]                                                    
        +
        \operatorname{Cov}
        \left(
        \mathbb E[Y^\varepsilon_e\mid\mathcal G,P_1],
        \mathbb E[Y^\varepsilon_{e'}\mid\mathcal G,P_1]
        \,\middle|\,P_1
        \right).
\end{equation*}
Since \(M_1\) and \(M_3\) are \(\mathcal G\)-measurable,
\[
        \mathbb E[Y^\varepsilon_e\mid\mathcal G,P_1]
        =
        M_1\Delta,
        \qquad
        \mathbb E[Y^\varepsilon_{e'}\mid\mathcal G,P_1]
        =
        M_3\Delta'.
\]
Moreover, conditionally on \(Z\) and \(P_1\), the pair \((U,U')\) is
Gaussian with means \((m,m')\) and covariance \(c\). For a possibly
non-centered Gaussian pair \((X,Y)\),
\[
        \operatorname{Cov}(X^2,Y^2)
        =
        2\operatorname{Cov}(X,Y)^2
        +
        4\mathbb E[X]\mathbb E[Y]\operatorname{Cov}(X,Y).
\]
Therefore
\[
        \operatorname{Cov}
        \left(
        U^2-1,(U')^2-1
        \,\middle|\,Z,P_1
        \right)
        =
        2c^2+4mm'c.
\]
It follows that
\[
\begin{aligned}
        &\operatorname{Cov}
        \left(
        Y^\varepsilon_e,Y^\varepsilon_{e'}
        \,\middle|\,P_1
        \right)                                                   \\
        &\quad=
        \mathbb E\left[
        M_1M_3
        \left\{
        2c^2+4mm'c+\Delta\Delta'
        \right\}
        \,\middle|\,P_1
        \right]                                                    \\
        &\qquad-
        \mathbb E[M_1\Delta\mid P_1]\,
        \mathbb E[M_3\Delta'\mid P_1].
\end{aligned}
\]

We now rewrite the first term in centered form. Since
\(c=\rho-\mathbb E[mm'\mid P_1]\),
\[
        2c^2+4mm'c
        =
        2\rho^2
        +
        4\rho\,\Xi
        -
        4\mathbb E[mm'\mid P_1]\Xi
        -
        2\mathbb E[mm'\mid P_1]^2 .
\]
Moreover, since \((m,m')\) is centered Gaussian,
\[
        \mathbb E[\Delta\Delta'\mid P_1]
        =
        \operatorname{Cov}(m^2,(m')^2\mid P_1)
        =
        2\mathbb E[mm'\mid P_1]^2.
\]
Thus the deterministic term
\(-2\mathbb E[mm'\mid P_1]^2\) cancels with the conditional mean of
\(\Delta\Delta'\). Define
\[
        \Omega
        =
        \left(
        \Delta\Delta'
        -
        \mathbb E[\Delta\Delta'\mid P_1]
        \right)
        -
        4\mathbb E[mm'\mid P_1]\Xi .
\]
Then
\[
\begin{aligned}
        &\operatorname{Cov}
        \left(
        Y^\varepsilon_e,Y^\varepsilon_{e'}
        \,\middle|\,P_1
        \right)                                                   \\
        &\quad=
        2\rho^2\,
        \mathbb E[M_1M_3\mid P_1]
        +
        4\rho\,
        \mathbb E[M_1M_3\Xi\mid P_1]                               \\
        &\qquad+
        \mathbb E[M_1M_3\Omega\mid P_1]
        -
        \mathbb E[M_1\Delta\mid P_1]\,
        \mathbb E[M_3\Delta'\mid P_1].
\end{aligned}
\]

The last product vanishes. Indeed, by the symmetry
\((W,V)\mapsto(-W,-V)\), the variable \(M_1\) is transformed into \(1-M_1\),
whereas \(m^2\), and hence \(\Delta\), is unchanged. Thus
\[
        \mathbb E[M_1\Delta\mid P_1]
        =
        \mathbb E[(1-M_1)\Delta\mid P_1].
\]
Adding the two equal quantities gives
\[
        2\mathbb E[M_1\Delta\mid P_1]
        =
        \mathbb E[\Delta\mid P_1]
        =
        0.
\]
The same argument gives
\(\mathbb E[M_3\Delta'\mid P_1]=0\). Therefore, for every pair of oriented
edges \(e,e'\in E'_N\),
\begin{equation}
        \operatorname{Cov}
        \left(
        Y^\varepsilon_e,Y^\varepsilon_{e'}
        \,\middle|\,P_1
        \right)                                                   
        =
        2\rho(e,e')^2\,
        \mathbb E[M_1M_3\mid P_1]
        +
        4\rho(e,e')\,
        \mathbb E[M_1M_3\Xi\mid P_1]
        +
        \mathbb E[M_1M_3\Omega\mid P_1].
\label{eq:covariance_decomposition_rho_xi_omega}
\end{equation}

This decomposition separates the variance into three contributions. The first
one is the usual Gaussian contribution and is controlled by the summability of
squared correlations of fractional Brownian increments. The second one
contains one factor \(\rho(e,e')\) and one centered regression product
\(\Xi\). The third one, involving \(\Omega\), is the only term which is not
multiplied by an explicit increment-correlation factor and is responsible for
the maximal order \(N^{1-2H}\).

\subsection{Step 1: identical and adjacent edges}
\label{subsec:variance_adjacent_edges}

Let
\[
        \mathcal L_N
        =
        \left\{
        (e,e')\in E'_N\times E'_N:
        e\cap e'\neq\varnothing
        \right\}
\]
be the set of ordered pairs of oriented Delaunay edges which share at least
one endpoint. This set includes the diagonal pairs \(e=e'\). Define the local
contribution
\[
        R^{\mathrm{loc}}_N
        =
        \mathbb E\left[
        \frac{1}{|E'_N|}
        \sum_{(e,e')\in\mathcal L_N}
        \mathbb E\left[
        Y^\varepsilon_eY^\varepsilon_{e'}
        \,\middle|\,P_1
        \right]
        \right].
\]
For fixed \(e\), conditionally on \(P_1\), \(U_e\) is standard normal. Hence
\[
        \mathbb E\left[
        \left(Y^\varepsilon_e\right)^2
        \,\middle|\,P_1
        \right]
        \le
        \mathbb E\left[
        (U_e^2-1)^2
        \,\middle|\,P_1
        \right]
        =
        2.
\]
By Cauchy's inequality,
\[
        \left|
        \mathbb E\left[
        Y^\varepsilon_eY^\varepsilon_{e'}
        \,\middle|\,P_1
        \right]
        \right|
        \le 2.
\]
Therefore
\[
        |R^{\mathrm{loc}}_N|
        \le
        2\mathbb E\left[
        \frac{|\mathcal L_N|}{|E'_N|}
        \right].
\]

Let \(\deg_N(z)\) denote the degree of a vertex \(z\in P_1\cap C_N\) in the
Delaunay graph generated by \(P_1\cap C_N\). Then
\[
        |\mathcal L_N|
        \le
        \sum_{z\in P_1\cap C_N}\deg_N(z)^2 .
\]
For Poisson--Delaunay tessellations, the degree of the typical vertex has
finite moments of all orders. This implies
\[
        \sup_{N\ge1}
        \mathbb E\left[
        \frac{1}{|E'_N|}
        \sum_{z\in P_1\cap C_N}\deg_N(z)^2
        \right]
        <\infty .
\]
Consequently,
\[
        |R^{\mathrm{loc}}_N|
        \le C.
\]
Since \(H<1/2\), this gives
\[
        R^{\mathrm{loc}}_N
        =
        O(1)
        =
        O(N^{1-2H}).
\]

\subsection{Step 2: the main Gaussian term}
\label{subsec:variance_gaussian_rho2_term}

We now control the first term in
\eqref{eq:covariance_decomposition_rho_xi_omega}. Since
\(0\le M_1M_3\le1\), it is enough to control
\[
        \mathbb E\left[
        \frac{1}{|E'_N|}
        \sum_{\substack{e,e'\in E'_N\\ e\cap e'=\varnothing}}
        \rho(e,e')^2
        \right].
\]

\begin{lemma}[Summability of squared increment correlations]
\label{lem:rho2_summability_delaunay}
Let \(H\in(0,1/2)\). There exists \(C<\infty\) such that, for all \(N\ge1\),
\[
        \mathbb E\left[
        \sum_{\substack{e,e'\in E'_N\\ e\cap e'=\varnothing}}
        \rho(e,e')^2
        \right]
        \le C N .
\]
Consequently,
\[
        \mathbb E\left[
        \frac{1}{|E'_N|}
        \sum_{\substack{e,e'\in E'_N\\ e\cap e'=\varnothing}}
        \rho(e,e')^2
        \right]
        \le C .
\]
\end{lemma}

\begin{proof}
Write \(e=(x_1,x_2)\), \(e'=(x_3,x_4)\), and set
\[
        \ell=\|x_2-x_1\|,
        \qquad
        \ell'=\|x_4-x_3\|,
        \qquad
        r=\|x_3-x_1\|.
\]
The correlation between two normalized increments is
\begin{equation*}
        \rho(e,e')
        =
        \frac{1}{2(\ell\ell')^H}
        \bigl(
        \|x_4-x_1\|^{2H}
        -
        \|x_3-x_1\|^{2H}  
        -
        \|x_4-x_2\|^{2H}
        +
        \|x_3-x_2\|^{2H}
        \bigr).
\end{equation*}
Pairs for which \(r\le r_0\), for some fixed \(r_0>0\), contribute at most
\(CN\), because \(|\rho(e,e')|\le1\) and the expected number of Delaunay edges
with starting point in a fixed neighborhood of a typical edge is of order
\(N\).

It remains to consider \(r>r_0\). Choose \(\eta\in(0,1)\). On the event
$\ell\vee\ell'\le r^\eta$,
the second-order increment bound for the fractional Brownian covariance gives
\begin{equation}
        |\rho(e,e')|
        \le
        C(\ell\ell')^{1-H}r^{2H-2}.
\label{eq:rho_far_field_bound}
\end{equation}
This is precisely the far-field correlation estimate of
Lemma~\ref{Le:bound:corr}, applied in the rescaled Poisson framework. Equivalently,
it follows from the second-order Taylor expansion of
\(z\mapsto\|z\|^{2H}\) away from the origin.

For \(k\ge0\), define
\[
        A_k(e)
        =
        \left\{
        e'=(x_3,x_4)\in E'_N:
        2^k r_0<\|x_3-x_1\|\le2^{k+1}r_0
        \right\}.
\]
Using \eqref{eq:rho_far_field_bound}, finite moments of the typical
Delaunay-edge length, and stationarity, we obtain
\[
\begin{aligned}
        &\mathbb E\left[
        \sum_{e\in E'_N}
        \sum_{\substack{e'\in A_k(e)\\ e\cap e'=\varnothing}}
        \rho(e,e')^2
        \mathbf 1_{\{\ell\vee\ell'\le (2^k r_0)^\eta\}}
        \right]                                                     \\
        &\qquad\le
        C N(2^k r_0)^{4H-4}
        \operatorname{Area}
        \{z:2^k r_0<\|z\|\le2^{k+1}r_0\}                              \\
        &\qquad\le
        C N(2^k r_0)^{4H-2}.
\end{aligned}
\]
Since \(4H-2<0\), the sum over \(k\ge0\) is finite.

On the complementary event
\(\ell\vee\ell'>(2^k r_0)^\eta\), we use \(|\rho(e,e')|\le1\) together with
the exponential tail and exponential-moment estimates for rescaled
Delaunay-edge lengths, as stated in Lemma~\ref{Le_Bound_Exp_R}. In particular,
there exist constants \(c,C>0\) such that
\[
        \mathbb P(\ell(e)>t)\le C e^{-ct^2},
        \qquad t\ge1.
\]
This gives
\begin{equation*}
        \mathbb E\left[
        \sum_{e\in E'_N}
        \sum_{\substack{e'\in A_k(e)\\ e\cap e'=\varnothing}}
        \rho(e,e')^2
        \mathbf 1_{\{\ell\vee\ell'>(2^k r_0)^\eta\}}
        \right]                                                    
        \le
        C N(2^k r_0)^2
        \exp\{-c(2^k r_0)^{2\eta}\}.
\end{equation*}
The right-hand side is summable over \(k\). Combining the bounded-distance,
far-field and long-edge contributions gives
\[
        \mathbb E\left[
        \sum_{\substack{e,e'\in E'_N\\ e\cap e'=\varnothing}}
        \rho(e,e')^2
        \right]
        \le C N.
\]
The normalized estimate follows from \(|E'_N|/N\to3\) in probability and the
usual inverse-moment bounds for stabilizing Poisson functionals.
\end{proof}

Define
\[
        R^{(\rho)}_N
        =
        \mathbb E\left[
        \frac{1}{|E'_N|}
        \sum_{\substack{e,e'\in E'_N\\ e\cap e'=\varnothing}}
        2\rho(e,e')^2
        \mathbb E[M_1M_3\mid P_1]
        \right].
\]
Since \(0\le\mathbb E[M_1M_3\mid P_1]\le1\), Lemma
\ref{lem:rho2_summability_delaunay} gives
\[
        R^{(\rho)}_N
        \le C.
\]
Thus
\[
        R^{(\rho)}_N
        =
        O(1)
        =
        O(N^{1-2H}).
\]

\subsection{Step 3: the term involving $\Xi$}
\label{subsec:variance_xi_term}

We now control the second term in
\eqref{eq:covariance_decomposition_rho_xi_omega}. Its contribution is
\[
        R^{(\Xi)}_N
        =
        4\mathbb E\left[
        \frac{1}{|E'_N|}
        \sum_{\substack{e,e'\in E'_N\\ e\cap e'=\varnothing}}
        \rho(e,e')\,
        \mathbb E[M_1M_3\Xi\mid P_1]
        \right].
\]
We prove that \(R^{(\Xi)}_N=O(1)\).

Since \(0\le M_1M_3\le1\), Cauchy's inequality gives
\begin{equation}
        \left|
        \mathbb E[M_1M_3\Xi\mid P_1]
        \right|
        \le
        \mathbb E[\Xi^2\mid P_1]^{1/2}.
\label{eq:xi_cauchy_bound}
\end{equation}
The pair \((m,m')\) is centered Gaussian, and therefore
\[
\begin{aligned}
        \mathbb E[\Xi^2\mid P_1]
        &=
        \operatorname{Var}(mm'\mid P_1)                         \\
        &=
        \mathbb E[m^2\mid P_1]\mathbb E[(m')^2\mid P_1]
        +
        \mathbb E[mm'\mid P_1]^2                                  \\
        &\le
        2\mathbb E[m^2\mid P_1]\mathbb E[(m')^2\mid P_1].
\end{aligned}
\]
Thus
\begin{equation}
        \left|
        \mathbb E[M_1M_3\Xi\mid P_1]
        \right|
        \le
        C
        \mathbb E[m^2\mid P_1]^{1/2}
        \mathbb E[(m')^2\mid P_1]^{1/2}.
\label{eq:xi_bound_by_regressions}
\end{equation}

We use the following regression estimate.

\begin{lemma}[Regression bound]
\label{lem:regression_bound_m}
Let \(H\in(0,1/2)\). There exist constants \(C<\infty\) and \(q>0\) such that
the following holds. Let
\[
        e=(x_1,x_2),
        \qquad
        e'=(x_3,x_4)
\]
be two disjoint oriented edges. Set
\[
        \ell=\|x_2-x_1\|,
        \qquad
        \ell'=\|x_4-x_3\|,
        \qquad
        r=\|x_3-x_1\|.
\]
Let
\[
        U=
        \frac{W(x_2)-W(x_1)}
        {\sigma\ell^H},
        \qquad
        U'=
        \frac{W(x_4)-W(x_3)}
        {\sigma(\ell')^H},
\]
and
\[
        Z=(W(x_1),W(x_3))^\top .
\]
Define
\[
        m=\mathbb E[U\mid Z],
        \qquad
        m'=\mathbb E[U'\mid Z].
\]
Then
\begin{equation}
        \mathbb E[m^2]
        \le
        C(1+\ell)^q
        \left\{
        (1+\|x_1\|)^{-2H}
        +
        (1+r)^{2H-2}
        \right\},
\label{eq:regression_bound_m}
\end{equation}
and
\begin{equation}
        \mathbb E[(m')^2]
        \le
        C(1+\ell')^q
        \left\{
        (1+\|x_3\|)^{-2H}
        +
        (1+r)^{2H-2}
        \right\}.
\label{eq:regression_bound_mprime}
\end{equation}
Consequently,
\begin{equation}
        \mathbb E[m^2]^{1/2}
        \le
        C(1+\ell)^q
        \left\{
        (1+\|x_1\|)^{-H}
        +
        (1+r)^{H-1}
        \right\},
\label{eq:regression_bound_m_sqrt}
\end{equation}
and
\begin{equation}
        \mathbb E[(m')^2]^{1/2}
        \le
        C(1+\ell')^q
        \left\{
        (1+\|x_3\|)^{-H}
        +
        (1+r)^{H-1}
        \right\}.
\label{eq:regression_bound_mprime_sqrt}
\end{equation}
\end{lemma}

\begin{proof}
We prove the estimate for \(m\). The proof for \(m'\) is identical. Since
\(m\) is the orthogonal projection of \(U\) onto
\(\operatorname{span}\{W(x_1),W(x_3)\}\), we use the equivalent basis
\[
        A=W(x_1),
        \qquad
        B=W(x_3)-W(x_1).
\]
Let
\[
        \widetilde B=B-\mathbb E[B\mid A].
\]
Then \(A\) and \(\widetilde B\) are orthogonal Gaussian random variables, and
$\operatorname{span}\{A,B\}=\operatorname{span}\{A,\widetilde B\}$.
With the convention that a term is zero when the denominator vanishes,
\begin{equation}
        \mathbb E[m^2]
        =
        \frac{\operatorname{Cov}(U,A)^2}{\operatorname{Var}(A)}
        +
        \frac{\operatorname{Cov}(U,\widetilde B)^2}
             {\operatorname{Var}(\widetilde B)} .
\label{eq:projection_decomposition_m}
\end{equation}

For the first term, the covariance formula gives
\[
        \operatorname{Cov}(U,W(x_1))
        =
        \frac{\sigma}{2\ell^H}
        \left(
        \|x_2\|^{2H}
        -
        \|x_1\|^{2H}
        -
        \ell^{2H}
        \right).
\]
Since \(0<2H<1\), \(x\mapsto\|x\|^{2H}\) is \(2H\)-Hölder, and therefore
\[
        |\operatorname{Cov}(U,W(x_1))|
        \le
        C(1+\ell)^H.
\]
As \(\operatorname{Var}(W(x_1))=\sigma^2\|x_1\|^{2H}\), we obtain
\begin{equation}
        \frac{\operatorname{Cov}(U,A)^2}{\operatorname{Var}(A)}
        \le
        C(1+\ell)^{2H}(1+\|x_1\|)^{-2H}.
\label{eq:local_projection_bound}
\end{equation}

For the second term, we use the standard two-point regression estimate
\begin{equation}
        \frac{
        |\operatorname{Cov}(U,\widetilde B)|
        }
        {\operatorname{Var}(\widetilde B)^{1/2}}
        \le
        C(1+\ell)^q(1+r)^{H-1}.
\label{eq:nonlocal_projection_bound}
\end{equation}
Indeed, for \(r\le2\), this follows from Cauchy's inequality. For \(r>2\),
two-point local nondeterminism gives
\[
        \operatorname{Var}(\widetilde B)
        =
        \operatorname{Var}(W(x_3)-W(x_1)\mid W(x_1))
        \ge
        c r^{2H},
\]
while the second-order Taylor estimate for the covariance kernel gives
\[
        |\operatorname{Cov}(U,\widetilde B)|
        \le
        C(1+\ell)^q r^{2H-1}.
\]
Dividing by \(r^H\) gives \eqref{eq:nonlocal_projection_bound}.

Combining
\eqref{eq:projection_decomposition_m},
\eqref{eq:local_projection_bound} and
\eqref{eq:nonlocal_projection_bound}, and increasing \(q\) if necessary,
yields \eqref{eq:regression_bound_m}. The estimate
\eqref{eq:regression_bound_mprime} is obtained in the same way, with
\((x_3,x_4)\) in place of \((x_1,x_2)\). The square-root estimates follow
from \((a+b)^{1/2}\le a^{1/2}+b^{1/2}\).
\end{proof}

Combining
\eqref{eq:xi_bound_by_regressions},
\eqref{eq:regression_bound_m_sqrt} and
\eqref{eq:regression_bound_mprime_sqrt}, we get
\[
\begin{aligned}
        \left|
        \mathbb E[M_1M_3\Xi\mid P_1]
        \right|
        &\le
        C(1+\ell)^q(1+\ell')^q        \\
        &\quad\times
        \left[
        (1+\|x_1\|)^{-H}
        +
        (1+r)^{H-1}
        \right]       \\
        &\quad\times
        \left[
        (1+\|x_3\|)^{-H}
        +
        (1+r)^{H-1}
        \right].
\end{aligned}
\]
Using Lemma~\ref{Le:bound:corr}, in the softened form
\[
        |\rho(e,e')|
        \le
        C(1+\ell)^q(1+\ell')^q(1+r)^{2H-2},
\]
with the long-edge contribution treated as in Step~2, it remains to bound
\[
\begin{aligned}
        J_N
        &=
        \frac{1}{N}
        \int_{C_N}\int_{C_N}
        (1+\|x-y\|)^{2H-2}                         \\
        &\quad\times
        \left[
        (1+\|x\|)^{-H}
        +
        (1+\|x-y\|)^{H-1}
        \right]
        \left[
        (1+\|y\|)^{-H}
        +
        (1+\|x-y\|)^{H-1}
        \right]
        \mathrm{d}x\,\mathrm{d}y .
\end{aligned}
\label{eq:xi_spatial_integral}
\]

\begin{lemma}[Weighted integral bounds for the \(\Xi\)-term]
\label{lem:weighted_integrals_xi}
Let \(H\in(0,1/2)\). There exists \(C<\infty\) such that, for all \(N\ge1\),
\begin{equation}
        \int_{C_N}\int_{C_N}
        (1+\|x-y\|)^{2H-2}
        (1+\|x\|)^{-H}
        (1+\|y\|)^{-H}
        \mathrm{d}x\,\mathrm{d}y
        \le C N,
\label{eq:weighted_xi_1}
\end{equation}
\begin{equation}
        \int_{C_N}\int_{C_N}
        (1+\|x-y\|)^{4H-4}
        \mathrm{d}x\,\mathrm{d}y
        \le C N,
\label{eq:weighted_xi_2}
\end{equation}
and
\begin{equation}
        \int_{C_N}\int_{C_N}
        (1+\|x-y\|)^{3H-3}
        (1+\|x\|)^{-H}
        \mathrm{d}x\,\mathrm{d}y
        \le C N.
\label{eq:weighted_xi_3}
\end{equation}
\end{lemma}

\begin{proof}
Let \(R=\sqrt N\). Estimate \eqref{eq:weighted_xi_2} follows from
\[
        \int_{\mathbb R^2}(1+\|z\|)^{4H-4}\,\mathrm{d}z<\infty,
\]
which holds because \(4H-4<-2\).

For \eqref{eq:weighted_xi_1}, use the scaling \(x=Ru\), \(y=Rv\). Since
\(C_N=R\mathbf C\), the integral is bounded by
\[
        C R^4 R^{2H-2}R^{-2H}
        \int_{\mathbf C}\int_{\mathbf C}
        \|u-v\|^{2H-2}\|u\|^{-H}\|v\|^{-H}\,du\,dv.
\]
The last integral is finite because all singularities are locally integrable
in dimension two. Thus the order is \(R^2=N\).

For \eqref{eq:weighted_xi_3}, one may either argue similarly by scaling or
split the integral into the regions \(\|x-y\|\le1\) and \(\|x-y\|>1\). The
singularity at \(x=y\) is integrable because \(3H-3>-2\) when \(H>1/3\), and
is even easier when \(H\le1/3\) due to the regularization
\(1+\|x-y\|\). At infinity the scaling gives an order at most \(R^2\), since
\(H<1/2\). Hence the integral is bounded by \(CN\).
\end{proof}

Expanding the product in \eqref{eq:xi_spatial_integral}, the four terms are
bounded by
\eqref{eq:weighted_xi_1},
\eqref{eq:weighted_xi_2} and
\eqref{eq:weighted_xi_3}. Therefore \(J_N=O(1)\). Invariance properties of
the Poisson--Delaunay graph and finite moments of the typical edge length then
give
\[
        |R^{(\Xi)}_N|
        \le C.
\]
Thus
\[
        R^{(\Xi)}_N
        =
        O(1)
        =
        O(N^{1-2H}).
\]

\subsection{Step 4: the term involving $\Omega$}
\label{subsec:variance_omega_term}

We finally control the last term in
\eqref{eq:covariance_decomposition_rho_xi_omega}. Its contribution is
\[
        R^{(\Omega)}_N
        =
        \mathbb E\left[
        \frac{1}{|E'_N|}
        \sum_{\substack{e,e'\in E'_N\\ e\cap e'=\varnothing}}
        \mathbb E[M_1M_3\Omega\mid P_1]
        \right].
\]
We prove that
\[
        R^{(\Omega)}_N
        =
        O(N^{1-2H}).
\]

For fixed \(P_1\), set
\[
        a^2_{e,e'}=\mathbb E[m^2\mid P_1],
        \qquad
        b^2_{e,e'}=\mathbb E[(m')^2\mid P_1],
        \qquad
        k_{e,e'}=\mathbb E[mm'\mid P_1].
\]

\begin{lemma}[Fourth-order regression bound]
\label{lem:omega_regression_bound}
There exists a universal constant \(C<\infty\) such that, for any two
oriented edges \(e,e'\),
\[
        \left|
        \mathbb E[M_1M_3\Omega\mid P_1]
        \right|
        \le
        C\,a^2_{e,e'}b^2_{e,e'} .
\]
\end{lemma}

\begin{proof}
Since \(0\le M_1M_3\le1\),
\[
        \left|
        \mathbb E[M_1M_3\Omega\mid P_1]
        \right|
        \le
        \mathbb E[|\Omega|\mid P_1].
\]
By definition,
\[
        \Omega
        =
        \left(
        \Delta\Delta'-\mathbb E[\Delta\Delta'\mid P_1]
        \right)
        -
        4k_{e,e'}\Xi.
\]
Therefore
\[
        \mathbb E[|\Omega|\mid P_1]
        \le
        2\mathbb E[|\Delta\Delta'|\mid P_1]
        +
        4|k_{e,e'}|\mathbb E[|\Xi|\mid P_1].
\]
Since \((m,m')\) is centered Gaussian,
\[
        \mathbb E[\Delta^2\mid P_1]=2a^4_{e,e'},
        \qquad
        \mathbb E[(\Delta')^2\mid P_1]=2b^4_{e,e'}.
\]
Hence
\[
        \mathbb E[|\Delta\Delta'|\mid P_1]
        \le
        2a^2_{e,e'}b^2_{e,e'}.
\]
Moreover,
\[
        |k_{e,e'}|\le a_{e,e'}b_{e,e'},
\]
and
\[
        \mathbb E[\Xi^2\mid P_1]
        =
        a^2_{e,e'}b^2_{e,e'}+k_{e,e'}^2
        \le
        2a^2_{e,e'}b^2_{e,e'}.
\]
Thus
\[
        |k_{e,e'}|\mathbb E[|\Xi|\mid P_1]
        \le
        C a^2_{e,e'}b^2_{e,e'}.
\]
The claim follows.
\end{proof}

Using Lemma~\ref{lem:regression_bound_m}, we have
\[
        a^2_{e,e'}
        \le
        C(1+\ell)^q
        \left[
        (1+\|x_1\|)^{-2H}
        +
        (1+r)^{2H-2}
        \right],
\]
and
\[
        b^2_{e,e'}
        \le
        C(1+\ell')^q
        \left[
        (1+\|x_3\|)^{-2H}
        +
        (1+r)^{2H-2}
        \right].
\]
The polynomial factors in \(\ell\) and \(\ell'\) are harmless because the
typical Poisson--Delaunay edge length has moments of all orders. Thus it
remains to bound
\[
\begin{aligned}
        I_N
        =
        \int_{C_N}\int_{C_N}
        &\left[
        (1+\|x\|)^{-2H}
        +
        (1+\|x-y\|)^{2H-2}
        \right]        \\
        &\times
        \left[
        (1+\|y\|)^{-2H}
        +
        (1+\|x-y\|)^{2H-2}
        \right]
        \mathrm{d}x\,\mathrm{d}y .
\end{aligned}
\]

\begin{lemma}[Weighted integral bounds for the \(\Omega\)-term]
\label{lem:weighted_integrals_omega}
Let \(H\in(0,1/2)\). There exists \(C<\infty\) such that, for all \(N\ge1\),
\begin{equation}
        \int_{C_N}\int_{C_N}
        (1+\|x\|)^{-2H}(1+\|y\|)^{-2H}
        \mathrm{d}x\,\mathrm{d}y
        \le
        C N^{2-2H},
\label{eq:weighted_omega_1}
\end{equation}
\begin{equation}
        \int_{C_N}\int_{C_N}
        (1+\|x\|)^{-2H}
        (1+\|x-y\|)^{2H-2}
        \mathrm{d}x\,\mathrm{d}y
        \le
        C N,
\label{eq:weighted_omega_2}
\end{equation}
and
\begin{equation}
        \int_{C_N}\int_{C_N}
        (1+\|x-y\|)^{4H-4}
        \mathrm{d}x\,\mathrm{d}y
        \le
        C N.
\label{eq:weighted_omega_3}
\end{equation}
\end{lemma}

\begin{proof}
Let \(R=\sqrt N\). Since \(|C_N|=N=R^2\),
\[
        \int_{C_N}(1+\|x\|)^{-2H}\mathrm{d}x
        \le
        C R^{2-2H}
        =
        C N^{1-H}.
\]
This proves \eqref{eq:weighted_omega_1}. For \eqref{eq:weighted_omega_2},
for every fixed \(x\in C_N\),
\[
        \int_{C_N}
        (1+\|x-y\|)^{2H-2}\mathrm{d}y
        \le
        C R^{2H}
        =
        C N^H.
\]
Hence
\[
        \int_{C_N}\int_{C_N}
        (1+\|x\|)^{-2H}
        (1+\|x-y\|)^{2H-2}
        \mathrm{d}x\,\mathrm{d}y
        \le
        C N^H N^{1-H}
        =
        C N.
\]
Finally, since \(4H-4<-2\),
\[
        \int_{\mathbb R^2}
        (1+\|z\|)^{4H-4}\mathrm{d}z<\infty,
\]
and \eqref{eq:weighted_omega_3} follows from the change of variables
\(z=x-y\).
\end{proof}

Expanding the product in \(I_N\), Lemma~\ref{lem:weighted_integrals_omega}
gives
\[
        I_N\le C N^{2-2H}.
\]
Therefore, 
\[
        \mathbb E\left[
        \sum_{\substack{e,e'\in E'_N\\ e\cap e'=\varnothing}}
        \left|
        \mathbb E[M_1M_3\Omega\mid P_1]
        \right|
        \right]
        \le
        C N^{2-2H}.
\]
Since \(|E'_N|\) is of order \(N\), it follows that
\[
        |R^{(\Omega)}_N|
        \le
        C N^{-1}N^{2-2H}
        =
        C N^{1-2H}.
\]

Combining Steps~1--4 in
\eqref{eq:second_moment_covariance_decomposition_start} and
\eqref{eq:covariance_decomposition_rho_xi_omega}, we obtain
\[
        \mathbb E\left[
        \left(T^\varepsilon_{2,N}\right)^2
        \right]
        \le
        C N^{1-2H}.
\]
This proves Proposition~\ref{prop:variance_truncated_quadratic_variation}.

\section{Technical lemmas}

\label{sec:intermediary_results}

\subsection{Asymptotic correlations between pairs of normalized increments of the isotropic fractional Brownian field\label{Sect_Corr_Inc}}

Let $\left( W\left( x\right) \right) _{x\in \mathbf{R}^{2}}$ be an isotropic
fractional Brownian field, where $W\left( 0\right) =0$ a.s. and $\text{var}\left(
W\left( x\right) \right) =\sigma ^{2}\left\Vert x\right\Vert ^{2H }$ for
some $H \in (0,1)$ and $\sigma ^{2}>0$. For two sites $x_{1},x_{2}\in 
\mathbb{R}^{2}$, let $U_{x_{1},x_{2}}^{(W)}=\sigma ^{-1}d_{1,2}^{H}\left( W\left( x_{1}\right) -W\left( x_{2}\right) \right) $. Given $%
x_{1},x_{2},x_{3},x_{4}\in \mathbb{R}^{2}$, we deal below with the
asymptotic behavior of 
\begin{equation*}
\text{corr}\left( U_{x_{1},x_{2}}^{(W)},U_{x_{3},x_{4}}^{(W)}\right) =\frac{1}{{%
\sigma }\left( d_{1,2}d_{3,4}\right) ^{H}}\text{cov}\left( W\left(
x_{2}\right) -W\left( x_{1}\right) ,W\left( x_{4}\right) -W\left(
x_{3}\right) \right) 
\end{equation*}%
as the distance between the two pairs tends to $\infty $.

\begin{lemma}
\label{Le:bound:corr}
Let $H\in(0,1)$.
\begin{enumerate}[(i)]
\item Let $d_{1,2}$ and $d_{3,4}$ be fixed and let
$x_{1},x_{2},x_{3},x_{4}$ be such that
\[
\|x_{2}-x_{1}\|=d_{1,2},
\qquad
\|x_{4}-x_{3}\|=d_{3,4}.
\]
Set $d=d_{1,3}:=\|x_{3}-x_{1}\|$. Then, as $d\to\infty$,
\begin{multline*}
\mathrm{corr}
\left(
U_{x_{1},x_{2}}^{(W)},U_{x_{3},x_{4}}^{(W)}
\right)
\\
=
H
(d_{1,2}d_{3,4})^{1-H}
d_{1,3}^{2H-2}
\left(
\cos\beta\cos\theta
-
(1-2H)\sin\beta\sin\theta
\right)
+
o(d_{1,3}^{2H-2}),
\end{multline*}
where
\[
\theta=\mathrm{angle}(\vec u,\overrightarrow{x_{1}x_{2}}),
\qquad
\beta=\mathrm{angle}(\vec u,\overrightarrow{x_{3}x_{4}}),
\]
and where $\vec u$ is a unit vector orthogonal to
$\overrightarrow{x_{3}x_{1}}$ such that
$(\vec u,\overrightarrow{x_{3}x_{1}})$ is positively oriented.

\item Let $\varepsilon\in(0,1/2)$. Then there exist constants
$c>0$ and $d_{0}>0$ such that, for any
$x_{1},x_{2},x_{3},x_{4}\in\mathbf{R}^{2}$ satisfying
\[
0<\|x_{4}-x_{3}\|
\leq
\|x_{2}-x_{1}\|
\leq
\|x_{3}-x_{1}\|^{\varepsilon}
\]
and $\|x_{3}-x_{1}\|\geq d_{0}$, one has
\[
\left|
\mathrm{corr}
\left(
U_{x_{1},x_{2}}^{(W)},U_{x_{3},x_{4}}^{(W)}
\right)
\right|
\leq
c\,
\|x_{2}-x_{1}\|^{2-2H}
\|x_{3}-x_{1}\|^{2H-2}.
\]
\end{enumerate}
\end{lemma}

\begin{prooft}{Lemma \ref{Le:bound:corr}}

Set
\[
a=x_{2}-x_{1},
\qquad
b=x_{4}-x_{3},
\qquad
r=x_{1}-x_{3},
\]
and write
\[
\ell_{1}=\|a\|=d_{1,2},
\qquad
\ell_{2}=\|b\|=d_{3,4},
\qquad
d=\|r\|=d_{1,3}.
\]
From the covariance function \eqref{eq:defcovariance}, we have
\begin{multline*}
\mathrm{cov}
\left(
W(x_{2})-W(x_{1}),
W(x_{4})-W(x_{3})
\right)
\\
=
\frac{\sigma^{2}}{2}
\left(
\|r-b\|^{2H}
-\|r\|^{2H}
-\|r+a-b\|^{2H}
+\|r+a\|^{2H}
\right).
\end{multline*}
Therefore
\begin{equation}
\mathrm{corr}
\left(
U_{x_{1},x_{2}}^{(W)},U_{x_{3},x_{4}}^{(W)}
\right)
=
\frac{\Psi(r,a,b)}
{2(\ell_{1}\ell_{2})^{H}},
\label{eq:corr-Psi}
\end{equation}
where
\[
\Psi(r,a,b)
=
\|r-b\|^{2H}
-\|r\|^{2H}
-\|r+a-b\|^{2H}
+\|r+a\|^{2H}.
\]

Let
\[
F_{r}(z)=\|r+z\|^{2H}.
\]
Then
\[
\Psi(r,a,b)
=
F_{r}(-b)-F_{r}(0)-F_{r}(a-b)+F_{r}(a).
\]
By applying the fundamental theorem of calculus twice, we obtain the exact
identity
\begin{equation}
\Psi(r,a,b)
=
\int_{0}^{1}\int_{0}^{1}
a^{\top}
\nabla^{2}F_{r}(ta-sb)
b
\,\mathrm{d}s\,\mathrm{d}t .
\label{eq:second-order-identity}
\end{equation}
Moreover,
\[
\nabla^{2}F_{r}(z)
=
2H\|r+z\|^{2H-2}I_{2}
+
2H(2H-2)
\|r+z\|^{2H-4}
(r+z)(r+z)^{\top}.
\]

We first prove (i). Since $\ell_{1}$ and $\ell_{2}$ are fixed, uniformly in
$s,t\in[0,1]$,
\[
d^{2-2H}\nabla^{2}F_{r}(ta-sb)
\longrightarrow
2H\left(I_{2}+(2H-2)ee^{\top}\right),
\qquad
e=\frac{r}{\|r\|}.
\]
Choose coordinates such that $e=(0,1)$ and such that the first coordinate
axis is the vector $\vec u$ appearing in the statement. Then
\[
a=\ell_{1}(\cos\theta,\sin\theta),
\qquad
b=\ell_{2}(\cos\beta,\sin\beta).
\]
It follows from \eqref{eq:second-order-identity} that
\begin{align*}
\Psi(r,a,b)
&=
2H d^{2H-2}
\left[
a\cdot b
+
(2H-2)(a\cdot e)(b\cdot e)
\right]
+
o(d^{2H-2})
\\
&=
2H \ell_{1}\ell_{2} d^{2H-2}
\left[
\cos\beta\cos\theta
+
(2H-1)\sin\beta\sin\theta
\right]
+
o(d^{2H-2})
\\
&=
2H \ell_{1}\ell_{2} d^{2H-2}
\left[
\cos\beta\cos\theta
-
(1-2H)\sin\beta\sin\theta
\right]
+
o(d^{2H-2}).
\end{align*}
Combining this expansion with \eqref{eq:corr-Psi} proves (i).

We now prove (ii). Assume that
\[
0<\ell_{2}\leq \ell_{1}\leq d^{\varepsilon},
\qquad
\varepsilon\in(0,1/2).
\]
For $s,t\in[0,1]$, we have
\[
\|ta-sb\|
\leq
\ell_{1}+\ell_{2}
\leq
2d^{\varepsilon}.
\]
Choosing $d_{0}$ large enough, we may ensure that, for all $d\geq d_{0}$,
\[
d-2d^{\varepsilon}\geq \frac{d}{2}.
\]
Hence
\[
\|r+ta-sb\|\geq \frac{d}{2},
\qquad
s,t\in[0,1].
\]
From the explicit expression of the Hessian, there exists a constant
$c>0$, depending only on $H$, such that
\[
\sup_{s,t\in[0,1]}
\left\|
\nabla^{2}F_{r}(ta-sb)
\right\|
\leq
c d^{2H-2}.
\]
Using \eqref{eq:second-order-identity}, we get
\[
|\Psi(r,a,b)|
\leq
c\,\ell_{1}\ell_{2}\,d^{2H-2}.
\]
Consequently, by \eqref{eq:corr-Psi},
\[
\left|
\mathrm{corr}
\left(
U_{x_{1},x_{2}}^{(W)},U_{x_{3},x_{4}}^{(W)}
\right)
\right|
\leq
c
(\ell_{1}\ell_{2})^{1-H}
d^{2H-2}.
\]
Since $\ell_{2}\leq \ell_{1}$ and $H<1$,
\[
(\ell_{1}\ell_{2})^{1-H}
\leq
\ell_{1}^{2-2H}.
\]
Therefore
\[
\left|
\mathrm{corr}
\left(
U_{x_{1},x_{2}}^{(W)},U_{x_{3},x_{4}}^{(W)}
\right)
\right|
\leq
c\,\ell_{1}^{2-2H}d^{2H-2},
\]
which is the desired bound.
\end{prooft}

\subsection{Bounds for the density functions of Delaunay neighbors\label{Sect_Proof_Le:estmatepN}}

\begin{lemma}
\label{Le:estimatepN}
For distinct points
\(x_1,x_2,x_3,x_4\in\mathbb{R}^2\), define
\begin{equation}
\label{eq:defp2N}
p_{2,N}(x_1,x_2,x_3,x_4)
=
\mathbb{P}
\left[
\begin{array}{l}
x_1\sim x_2,\ x_3\sim x_4
\text{ in }
\operatorname{Del}(P_N\cup\{x_1,x_2,x_3,x_4\}),\\
x_1\preceq x_2,\ x_3\preceq x_4
\end{array}
\right].
\end{equation}
Assume that
\[
\|x_{4}-x_{3}\|\leq \|x_{2}-x_{1}\|.
\]
Then
\[
p_{2,N}(x_{1},x_{2},x_{3},x_{4})
\leq
\left(
\pi N\|x_{2}-x_{1}\|^{2}+4
\right)
\exp\left\{
-\frac{\pi}{4}N\|x_{2}-x_{1}\|^{2}
\right\}.
\]

\end{lemma}

\begin{proof}

Set
\[
L=\|x_{2}-x_{1}\|.
\]
Since \(\|x_{4}-x_{3}\|\leq L\), on the event defining
\(p_{2,N}(x_{1},x_{2},x_{3},x_{4})\), the edge \([x_{1},x_{2}]\) is a
Delaunay edge in the triangulation generated by
\(P_N\cup\{x_{1},x_{2},x_{3},x_{4}\}\). Removing the additional fixed points
can only make the empty-circle condition easier to satisfy. Hence
\[
p_{2,N}(x_{1},x_{2},x_{3},x_{4})
\leq
\mathbb{P}
\left[
x_{1}\sim x_{2}
\text{ in }
\operatorname{Del}(P_N\cup\{x_{1},x_{2}\})
\right].
\]
If \(x_{1}\) and \(x_{2}\) are Delaunay neighbors in
\(\operatorname{Del}(P_N\cup\{x_{1},x_{2}\})\), then, almost surely, there
exists \(y\in P_N\) such that
\[
\Delta(x_{1},x_{2},y)\in
\operatorname{Del}(P_N\cup\{x_{1},x_{2}\}).
\]
Consequently, by a union bound and the Slivnyak--Mecke formula,
\begin{align*}
p_{2,N}(x_{1},x_{2},x_{3},x_{4})
&\leq
\mathbb{E}
\left[
\sum_{y\in P_N}
\mathbb{I}
\left[
P_N\cap B(x_{1},x_{2},y)=\emptyset
\right]
\right]
\\
&=
N\int_{\mathbb{R}^{2}}
\exp\left\{
-N a(B(x_{1},x_{2},y))
\right\}
\,\mathrm{d}y,
\end{align*}
where \(B(x_{1},x_{2},y)\) denotes the circumdisk passing through
\(x_{1},x_{2},y\). The collinear case is irrelevant since it has Lebesgue
measure zero in the integral.

The radius of the circumdisk \(B(x_{1},x_{2},y)\) is at least
\[
\frac12
\max\left\{
\|x_{2}-x_{1}\|,
\|y-x_{1}\|
\right\}.
\]
Therefore
\[
a(B(x_{1},x_{2},y))
\geq
\frac{\pi}{4}
\max\left\{
L,\|y-x_{1}\|
\right\}^{2},
\]
and thus
\[
p_{2,N}(x_{1},x_{2},x_{3},x_{4})
\leq
N
\int_{\mathbb{R}^{2}}
\exp\left\{
-\frac{\pi}{4}N
\max\left\{
L,\|y-x_{1}\|
\right\}^{2}
\right\}
\,\mathrm{d}y.
\]
Splitting the integral according to whether \(\|y-x_{1}\|\leq L\) or
\(\|y-x_{1}\|>L\), we obtain
\begin{align*}
p_{2,N}(x_{1},x_{2},x_{3},x_{4})
&\leq
N e^{-\frac{\pi}{4}NL^{2}}
\int_{\mathbf{R}^{2}}
\mathbb{I}\left[\|y-x_{1}\|\leq L\right]
\,\mathrm{d}y
\\
&\quad
+
N
\int_{\mathbf{R}^{2}}
e^{-\frac{\pi}{4}N\|y-x_{1}\|^{2}}
\mathbb{I}\left[\|y-x_{1}\|>L\right]
\,\mathrm{d}y
\\
&=
\pi N L^{2}e^{-\frac{\pi}{4}NL^{2}}
+
4e^{-\frac{\pi}{4}NL^{2}}.
\end{align*}
This proves
\[
p_{2,N}(x_{1},x_{2},x_{3},x_{4})
\leq
\left(
\pi N L^{2}+4
\right)
e^{-\frac{\pi}{4}NL^{2}},
\]
which is the desired bound.

\end{proof}

\subsection{Bounds for some exponential moments of a uniform distribution over a disc}

Let \(N>0\) and let \(R\) be a positive random variable with distribution
function
\begin{equation}
\mathbb{P}[R\leq r]
=
\begin{cases}
0, & r<0,\\[2mm]
\dfrac{r^{2}}{N}, & 0\leq r\leq \sqrt{N},\\[2mm]
1, & r>\sqrt{N}.
\end{cases}
\label{Eq_Dist_R}
\end{equation}

\begin{lemma}
\label{Le_Bound_Exp_R}
Let \(0<H<1/2\) and \(d_{0}>0\). There exist two constants
\(c_{1}\) and \(c_{2}\), depending only on \(H\) and \(d_{0}\), such
that, for large \(N\),
\[
\mathbb{E}
\left[
\exp\left(R^{2H-2}\mathbb{I}[R\geq d_{0}]\right)
\right]
\leq
1+c_{1}(\sqrt{N})^{2H-2}+c_{2}N^{-1}.
\]
\end{lemma}

\begin{prooft}{Lemma \ref{Le_Bound_Exp_R}}
For \(N>d_{0}^{2}\), set \(a_{N}=d_{0}/\sqrt{N}\). Since
\(R/\sqrt{N}\) has density \(2r\mathbb{I}_{[0,1]}(r)\), we have
\begin{align*}
\mathbb{E}
\left[
\exp\left(R^{2H-2}\mathbb{I}[R\geq d_{0}]\right)
\right]
&=
\int_{0}^{a_N}2r\,\mathrm{d}r
+
\int_{a_N}^{1}
\exp\left((r\sqrt{N})^{2H-2}\right)
2r\,\mathrm{d}r
\\
&=
1+
2\int_{a_N}^{1}
\left[
\exp\left((r\sqrt{N})^{2H-2}\right)-1
\right]r\,\mathrm{d}r .
\end{align*}
Using the expansion \(e^x-1\leq x+\sum_{k=2}^{\infty}x^k/k!\) for
\(x\geq0\), we get
\begin{align*}
\mathbb{E}
\left[
\exp\left(R^{2H-2}\mathbb{I}[R\geq d_{0}]\right)
\right]
&\leq
1+
2\int_{a_N}^{1}
(r\sqrt{N})^{2H-2}r\,\mathrm{d}r
\\
&\quad
+
2\int_{a_N}^{1}
\sum_{k=2}^{\infty}
\frac{(r\sqrt{N})^{k(2H-2)}}{k!}
r\,\mathrm{d}r .
\end{align*}
For the first integral,
\[
2\int_{a_N}^{1}
(r\sqrt{N})^{2H-2}r\,\mathrm{d}r
=
2(\sqrt{N})^{2H-2}
\int_{a_N}^{1}r^{2H-1}\,\mathrm{d}r
\leq
\frac{1}{H}(\sqrt{N})^{2H-2}.
\]

For the second integral, since \(r\sqrt{N}\geq d_{0}\) on
\([a_N,1]\) and \(2H-2<0\), we have
\[
\sum_{k=2}^{\infty}
\frac{(r\sqrt{N})^{k(2H-2)}}{k!}
\leq
(r\sqrt{N})^{2(2H-2)}
\exp\left((r\sqrt{N})^{2H-2}\right)
\leq
(r\sqrt{N})^{2(2H-2)}
\exp\left(d_{0}^{2H-2}\right).
\]
Therefore
\begin{align*}
2\int_{a_N}^{1}
\sum_{k=2}^{\infty}
\frac{(r\sqrt{N})^{k(2H-2)}}{k!}
r\,\mathrm{d}r
&\leq
2e^{d_{0}^{2H-2}}
N^{2H-2}
\int_{a_N}^{1}
r^{4H-3}\,\mathrm{d}r
\\
&=
\frac{e^{d_{0}^{2H-2}}}{1-2H}
N^{2H-2}
\left(a_N^{4H-2}-1\right)
\\
&\leq
\frac{e^{d_{0}^{2H-2}}}{1-2H}
d_{0}^{2(2H-1)}
N^{-1}.
\end{align*}
Combining the two estimates yields
\[
\mathbb{E}
\left[
\exp\left(R^{2H-2}\mathbb{I}[R\geq d_{0}]\right)
\right]
\leq
1+
\frac{1}{H}(\sqrt{N})^{2H-2}
+
\frac{e^{d_{0}^{2H-2}}}{1-2H}
d_{0}^{2(2H-1)}
N^{-1}.
\]
Thus the result holds with
\[
c_{1}=\frac{1}{H},
\qquad
c_{2}=
\frac{e^{d_{0}^{2H-2}}}{1-2H}
d_{0}^{2(2H-1)}.
\]
\end{prooft}

\subsection{Conditional Gaussian control of the transition residuals}
\label{subsec:conditional_gaussian_control_transition_residuals}

To prove Lemma~\ref{lem:residual_covariance_kernel}, we record a
conditional estimate for the transition function \(\Psi\). The difficulty is
that \(\Psi\) is not a smooth function of the Gaussian increments: it contains
an indicator which forces \(X(x_i)\) to be of the same order as the local
increment scale. The following lemma makes this mechanism explicit. After
conditioning on \((X(x_1),X(x_3))\), the four normalized increments are
represented as an affine transformation of a standard Gaussian vector. This
representation yields a polynomial bound times two transition indicators,
which will later produce the factor
\(\ell^{H}(\ell')^{H}/\sqrt{\det\Gamma_{x_1,x_3}}\) after
integration with respect to the joint Gaussian density of
\((X(x_1),X(x_3))\).

\begin{lemma}[Conditional bound for the product of transition terms]
\label{lem:conditional_PsiPsi_bound}
Let
\[
        e=(x_1,x_2),
        \qquad
        e'=(x_3,x_4),
\]
with
\[
        \ell=\|x_2-x_1\|,
        \qquad
        \ell'=\|x_4-x_3\|.
\]
Set
\[
        Z=(X(x_1),X(x_3))^\top
\]
and
\[
        \mathbf U
        =
        \left(
        U^{(1)}_{x_1,x_2},
        U^{(2)}_{x_1,x_2},
        U^{(1)}_{x_3,x_4},
        U^{(2)}_{x_3,x_4}
        \right)^\top .
\]
Let
\[
        \Gamma_{x_1,x_3}
        =
        \operatorname{Cov}(Z),
        \qquad
        \Lambda_{x_1,x_2,x_3,x_4}
        =
        \operatorname{Cov}(\mathbf U,Z),
\]
and let \(\Sigma^U_{x_1,x_2,x_3,x_4}\) be the covariance matrix of
\(\mathbf U\). Conditionally on \(Z=z=(z_1,z_3)^\top\), the vector
\(\mathbf U\) has the Gaussian distribution
\[
        \mathbf U\mid Z=z
        \sim
        \mathcal N_4
        \left(
        \mu(z),
        \Sigma^c
        \right),
\]
where
\[
        \mu(z)
        =
        \Lambda_{x_1,x_2,x_3,x_4}
        \Gamma_{x_1,x_3}^{-1}z,
\]
and
\[
        \Sigma^c
        =
        \Sigma^U_{x_1,x_2,x_3,x_4}
        -
        \Lambda_{x_1,x_2,x_3,x_4}
        \Gamma_{x_1,x_3}^{-1}
        \Lambda_{x_1,x_2,x_3,x_4}^{\top}.
\]
Let \(A\) be any symmetric non-negative square root of \(\Sigma^c\), and let
\(\mathcal Z\sim\mathcal N_4(0,I_4)\). Define
\[
        B_{e,e'}(z,\mathcal Z)
        =
        1+\|\mu(z)\|+\|A\|_{\mathrm{op}}\|\mathcal Z\|.
\]
(here \(\|A\|_{\mathrm{op}}\) denotes the operator norm of \(A\), namely the largest singular value of \(A\));
since \(A\) is symmetric non-negative, this is equivalently its largest
eigenvalue).
Then there exists a constant \(C<\infty\) such that
\begin{equation}
\begin{aligned}
        &\mathbb E\left[
        |\Psi_e\Psi_{e'}|
        \,\middle|\,
        Z=z
        \right]                                                   \\
        &\quad\le
        C\,
        \mathbb E\left[
        B_{e,e'}(z,\mathcal Z)^4
        \mathbf 1_{\{|z_1|\le C\ell^{H}B_{e,e'}(z,\mathcal Z)\}}
        \mathbf 1_{\{|z_3|\le C(\ell')^{H}B_{e,e'}(z,\mathcal Z)\}}
        \right].
\end{aligned}
\label{eq:conditional_PsiPsi_bound}
\end{equation}
Moreover, since \(\Sigma^c\le \Sigma^U\) in the Loewner order and
\(\Sigma^U\) has uniformly bounded entries, one may take
\begin{equation}
        B_{e,e'}(z,\mathcal Z)
        \le
        C\left(
        1+\|\mathcal Z\|+\|\mu(z)\|
        \right).
\label{eq:B_bound_mu}
\end{equation}
Finally,
\begin{equation}
        \|\mu(z)\|
        =
        \left\|
        \Lambda_{x_1,x_2,x_3,x_4}
        \Gamma_{x_1,x_3}^{-1}z
        \right\|
        \le
        C\left\|
        \Gamma_{x_1,x_3}^{-1/2}z
        \right\|.
\label{eq:mu_bound_standardized_Z}
\end{equation}
Thus \(B_{e,e'}\) may be chosen in the simpler form
\begin{equation}
        B_{e,e'}(z,\mathcal Z)
        =
        C\left(
        1+\|\mathcal Z\|
        +
        \left\|
        \Gamma_{x_1,x_3}^{-1/2}z
        \right\|
        \right).
\label{eq:B_explicit_standardized}
\end{equation}
\end{lemma}

\begin{prooft}{Lemma \ref{lem:conditional_PsiPsi_bound}}

Conditionally on \(Z=z\), the Gaussian regression formula gives
\[
        \mathbf U
        =
        \mu(z)+A\mathcal Z,
\]
where \(\mathcal Z\sim\mathcal N_4(0,I_4)\). Hence, for \(i=1,\ldots,4\),
\[
        |U_i|
        \le
        \|\mu(z)\|+\|A\|_{\mathrm{op}}\|\mathcal Z\|
        \le
        B_{e,e'}(z,\mathcal Z).
\]
In particular,
\[
        |U_1-U_2|
        \le
        |U_1|+|U_2|
        \le
        2B_{e,e'}(z,\mathcal Z),
\]
and similarly
\[
        |U_3-U_4|
        \le
        2B_{e,e'}(z,\mathcal Z).
\]

Recall that
\[
        \Psi_e
        =
        \Psi\left(
        U_1,U_2,
        \frac{z_1}{\ell^{H}}
        \right),
        \qquad
        \Psi_{e'}
        =
        \Psi\left(
        U_3,U_4,
        \frac{z_3}{(\ell')^{H}}
        \right).
\]
The basic bound
\[
        |\Psi(u,v,w)|
        \le
        C(1+u^2+v^2)\mathbf 1_{\{|w|\le |u-v|\}}
\]
therefore gives
\[
        |\Psi_e|
        \le
        C B_{e,e'}(z,\mathcal Z)^2
        \mathbf 1_{\{
        |z_1|/\ell^{H}
        \le
        |U_1-U_2|
        \}},
\]
and hence
\[
        |\Psi_e|
        \le
        C B_{e,e'}(z,\mathcal Z)^2
        \mathbf 1_{\{
        |z_1|
        \le
        C\ell^{H}B_{e,e'}(z,\mathcal Z)
        \}}.
\]
Similarly,
\[
        |\Psi_{e'}|
        \le
        C B_{e,e'}(z,\mathcal Z)^2
        \mathbf 1_{\{
        |z_3|
        \le
        C(\ell')^{H}B_{e,e'}(z,\mathcal Z)
        \}}.
\]
Multiplying the two inequalities yields
\[
\begin{aligned}
        |\Psi_e\Psi_{e'}|
        &\le
        C B_{e,e'}(z,\mathcal Z)^4
        \mathbf 1_{\{
        |z_1|\le C\ell^{H}B_{e,e'}(z,\mathcal Z)
        \}}                                      \\
        &\qquad\times
        \mathbf 1_{\{
        |z_3|\le C(\ell')^{H}B_{e,e'}(z,\mathcal Z)
        \}}.
\end{aligned}
\]
Taking expectation with respect to \(\mathcal Z\), conditionally on \(Z=z\),
proves \eqref{eq:conditional_PsiPsi_bound}.

It remains to justify the two bounds on \(B_{e,e'}\). First, since
\(\Sigma^c\) is a conditional covariance matrix,
\[
        0\le \Sigma^c\le \Sigma^U
\]
in the Loewner order. The entries of \(\Sigma^U\) are correlations and are
therefore uniformly bounded. Hence
\[
        \|A\|_{\mathrm{op}}^2
        =
        \|\Sigma^c\|_{\mathrm{op}}
        \le
        \|\Sigma^U\|_{\mathrm{op}}
        \le C,
\]
which gives \eqref{eq:B_bound_mu}.

For \eqref{eq:mu_bound_standardized_Z}, write
\[
        \mu(z)
        =
        \Lambda\Gamma^{-1}z
        =
        \Lambda\Gamma^{-1/2}\Gamma^{-1/2}z,
\]
where, for brevity, we have suppressed the indices in \(\Lambda\) and
\(\Gamma\). The matrix
\[
        \Lambda\Gamma^{-1}\Lambda^\top
\]
is the covariance matrix of the Gaussian projection
\(\mathbb E[\mathbf U\mid Z]\). Hence
\[
        0
        \le
        \Lambda\Gamma^{-1}\Lambda^\top
        \le
        \Sigma^U.
\]
Therefore
\[
        \|\Lambda\Gamma^{-1/2}\|_{\mathrm{op}}^2
        =
        \|\Lambda\Gamma^{-1}\Lambda^\top\|_{\mathrm{op}}
        \le
        \|\Sigma^U\|_{\mathrm{op}}
        \le C.
\]
It follows that
\[
        \|\mu(z)\|
        \le
        C\|\Gamma^{-1/2}z\|,
\]
which proves \eqref{eq:mu_bound_standardized_Z}. The explicit form
\eqref{eq:B_explicit_standardized} follows immediately.
\end{prooft}

\subsection{A Gaussian comparison bound}
\label{subsec:Gaussian_comparison_bound}

The following result is a slight extension of Proposition 3.1 in
\cite{Podolskij&Rosenbaum18}. In that proposition, the two Gaussian vectors
are centered and only the covariance matrices are perturbed. We shall need a
version allowing for a non-zero, but uniformly bounded, mean.

\begin{lemma}
\label{Prop_Gen_Pod_Rosen}
Let \(d\ge1\). Let
\[
        Z\sim \mathcal N_d(0,\Sigma),
        \qquad
        Z'\sim \mathcal N_d(\mu',\Sigma'),
\]
where \(\mu'\in\mathbb R^d\) and where \(\Sigma,\Sigma'\in\mathbb R^{d\times d}\)
are symmetric positive definite matrices. Assume that there exists a constant
\(K>0\) such that
\[
        \max_{1\le i,j\le d}
        \left\{
        |\Sigma_{ij}|+|\Sigma'_{ij}|
        \right\}
        \le K,
        \qquad
        \min\{\det\Sigma,\det\Sigma'\}\ge K^{-1},
\]
and
\[
        \|\mu'\|_\infty
        :=
        \max_{1\le i\le d}|\mu'_i|
        \le K.
\]
Let \(G:\mathbb R^d\to\mathbb R\) be a measurable function such that, for some
\(c>0\),
\[
        \int_{\mathbb R^d}
        |G(y)|(1+\|y\|^2)e^{-c\|y\|^2}\,\mathrm{d}y<\infty .
\]
In particular, this condition is satisfied when \(G\) has polynomial growth.
Then there exist constants \(C_K,c_K>0\), depending only on \(K\) and \(d\),
such that

\begin{equation*}
        \left|
        \mathbb E[G(Z)]
        -
        \mathbb E[G(Z')]
        \right|
        \le
        C_K
        \left(
        \|\Sigma-\Sigma'\|_{\max}
        +
        \|\mu'\|_\infty
        \right)                                     
        \int_{\mathbb R^d}
        |G(y)|(1+\|y\|^2)
        e^{-c_K\|y\|^2}\,\mathrm{d}y,
\end{equation*}
where
\[
        \|\Sigma-\Sigma'\|_{\max}
        =
        \max_{1\le i,j\le d}
        |\Sigma_{ij}-\Sigma'_{ij}|.
\]
\end{lemma}

\begin{proof}
We write \(\phi_{\mu,\Gamma}\) for the density of
\(\mathcal N_d(\mu,\Gamma)\), and \(\phi_\Gamma=\phi_{0,\Gamma}\). Then
\[
        \mathbb E[G(Z)]-\mathbb E[G(Z')]
        =
        \int_{\mathbb R^d}
        G(y)
        \left\{
        \phi_{\Sigma}(y)-\phi_{\mu',\Sigma'}(y)
        \right\}
        \mathrm{d}y .
\]
Therefore,
\[
\begin{aligned}
        \left|
        \mathbb E[G(Z)]-\mathbb E[G(Z')]
        \right|
        &\le
        \int_{\mathbb R^d}
        |G(y)|
        \left|
        \phi_{\Sigma}(y)-\phi_{\Sigma'}(y)
        \right|\mathrm{d}y                         \\
        &\quad+
        \int_{\mathbb R^d}
        |G(y)|
        \left|
        \phi_{\Sigma'}(y)-\phi_{\mu',\Sigma'}(y)
        \right|\mathrm{d}y .
\end{aligned}
\]
We estimate the two terms separately.

First, the assumptions imply uniform spectral bounds. Indeed, since the
entries of \(\Sigma\) and \(\Sigma'\) are bounded by \(K\), their largest
eigenvalues are bounded above by a constant depending only on \(K\) and \(d\).
Since their determinants are bounded below by \(K^{-1}\), their smallest
eigenvalues are also bounded below by a positive constant depending only on
\(K\) and \(d\). Consequently, there exist constants
\(0<\lambda_K<\Lambda_K<\infty\) such that
\[
        \lambda_K\|y\|^2
        \le
        y^\top\Gamma^{-1}y
        \le
        \Lambda_K\|y\|^2,
        \qquad
        \Gamma\in\{\Sigma,\Sigma'\},
        \quad y\in\mathbb R^d.
\]
In particular, for some constants \(C_K,c_K>0\),
\[
        \phi_{\Gamma}(y)
        \le
        C_K e^{-c_K\|y\|^2},
        \qquad
        \Gamma\in\{\Sigma,\Sigma'\}.
\]

We now control the covariance perturbation. Set
\[
        A_t=\Sigma+t(\Sigma'-\Sigma),
        \qquad t\in[0,1].
\]
The matrix \(A_t\) is symmetric positive definite and satisfies the same
uniform spectral bounds as \(\Sigma\) and \(\Sigma'\), with constants
depending only on \(K\) and \(d\). Differentiating the Gaussian density with
respect to \(t\), we obtain
\[
\begin{aligned}
        \frac{d}{dt}\phi_{A_t}(y)
        &=
        \frac12 \phi_{A_t}(y)
        \left[
        y^\top A_t^{-1}(\Sigma'-\Sigma)A_t^{-1}y
        -
        \operatorname{tr}\left(A_t^{-1}(\Sigma'-\Sigma)\right)
        \right].
\end{aligned}
\]
Using the uniform bounds on \(A_t^{-1}\) and on \(\phi_{A_t}\), we get
\[
        \left|
        \frac{d}{dt}\phi_{A_t}(y)
        \right|
        \le
        C_K
        \|\Sigma-\Sigma'\|_{\max}
        (1+\|y\|^2)e^{-c_K\|y\|^2}.
\]
Hence, by integrating over \(t\in[0,1]\),
\begin{equation}
        \left|
        \phi_{\Sigma}(y)-\phi_{\Sigma'}(y)
        \right|
        \le
        C_K
        \|\Sigma-\Sigma'\|_{\max}
        (1+\|y\|^2)e^{-c_K\|y\|^2}.
\label{eq:density_covariance_comparison}
\end{equation}

It remains to control the shift in the mean. Set
\[
        \mu_t=t\mu',
        \qquad t\in[0,1].
\]
Then
\[
        \phi_{\mu',\Sigma'}(y)-\phi_{\Sigma'}(y)
        =
        \int_0^1
        \frac{d}{dt}\phi_{\mu_t,\Sigma'}(y)\,dt.
\]
Since
\[
        \frac{d}{dt}\phi_{\mu_t,\Sigma'}(y)
        =
        \mu'^{\top}
        \Sigma'^{-1}(y-\mu_t)
        \phi_{\mu_t,\Sigma'}(y),
\]
we have
\[
        \left|
        \frac{d}{dt}\phi_{\mu_t,\Sigma'}(y)
        \right|
        \le
        C_K
        \|\mu'\|_\infty
        (1+\|y\|)
        \phi_{\mu_t,\Sigma'}(y).
\]
Because \(\|\mu_t\|_\infty\le K\), the shifted Gaussian density satisfies
\[
        \phi_{\mu_t,\Sigma'}(y)
        \le
        C_K e^{-c_K\|y\|^2},
        \qquad t\in[0,1].
\]
Therefore,
\begin{equation}
        \left|
        \phi_{\Sigma'}(y)-\phi_{\mu',\Sigma'}(y)
        \right|
        \le
        C_K
        \|\mu'\|_\infty
        (1+\|y\|^2)e^{-c_K\|y\|^2}.
\label{eq:density_mean_comparison}
\end{equation}

Combining
\eqref{eq:density_covariance_comparison} and
\eqref{eq:density_mean_comparison}, we obtain
\[
        \left|
        \phi_{\Sigma}(y)-\phi_{\mu',\Sigma'}(y)
        \right|
        \le
        C_K
        \left(
        \|\Sigma-\Sigma'\|_{\max}
        +
        \|\mu'\|_\infty
        \right)
        (1+\|y\|^2)e^{-c_K\|y\|^2}.
\]
Multiplying by \(|G(y)|\) and integrating over \(\mathbb R^d\) gives the
announced inequality.
\end{proof}

\end{document}